%% file: dichotomyManuscript.tex
\documentclass[preprint,12pt]{elsarticle}
 
\usepackage[margin=1in]{geometry} 
\usepackage{amsmath,amsthm,amssymb,scrextend}
\usepackage{mathtools}
\usepackage{fancyhdr}
\input{commands.tex}

\usepackage{tikz}
\usepackage{pgfplots}
\usepackage{amsmath}
\usepackage[mathscr]{euscript}
 \let\mathscr\relax% just so we can load this and rsfs
\usepackage[scr]{rsfso}
\usepackage{amsthm}
\usepackage[title]{appendix}
\usepackage{amssymb}
\usepackage{multicol}
\usepackage{soul}
\usepackage{changepage}
\usepackage[colorlinks=true, pdfstartview=FitV, linkcolor=blue,
citecolor=blue, urlcolor=blue]{hyperref}

\newcommand{\Ex}{\mathbb{E}}
\newcommand{\pr}[1]{\Pr\left[#1\right]}
\newcommand{\Var}{\operatorname{Var}}
\newcommand{\floor}[1]{\left\lfloor {#1} \right\rfloor}
\newcommand{\ceil}[1]{\left\lceil {#1} \right\rceil}

\newtheorem{theorem}{Theorem}
\numberwithin{theorem}{section}
\newtheorem{definition}[theorem]{Definition}
\newtheorem{lemma}[theorem]{Lemma}
\newtheorem{observation}[theorem]{Observation}

\newtheorem{corollary}[theorem]{Corollary}

\newtheorem{remk}[theorem]{Remark}

\newcommand{\N}{{\mathbb{N}}}

\newcommand{\dm}{d_{min}}
\newcommand{\f}{f(n,k, \dmax)}
\newcommand{\fopt}{f^*(n, \dmax)}

\newcommand{\gopt}{g^*(n, \dmax)}
\begin{document}

\begin{frontmatter}



\title{A dichotomy theorem on the complexity of 3-uniform hypergraphic degree sequence graphicality}


%Remark (Istvan): usual order of authors is alphabetic order.

\author[sara]{Sara Logsdon} %% Author name

\author[arya]{Arya Maheshwari}

\author[renyi,sztaki]{Istv\'an Mikl\'os}

\author[angelina]{Angelina Zhang}


\affiliation[sara]{organization={Department of Mathematics, University of Georgia},%Department and Organization
            addressline={200 D. W. Brooks Drive}, 
            city={Athens}, 
            state={GA},
            postcode={30602},
            country={USA}}

\affiliation[arya]{organization={Department of Computer Science, Princeton University},%Department and Organization
            addressline={35 Olden St}, 
            city={Princeton}, 
            state={NJ},
            postcode={08540},
            country={USA}}
\affiliation[renyi]{organization={HUN-REN Rényi Institute},%Department and Organization
            addressline={Reáltanoda u. 13--15}, 
            city={Budapest},
            postcode={1053}, 
            country={Hungary}}
\affiliation[sztaki]{organization={HUN-REN SZTAKI},%Department and Organization
            addressline={Lágymányosi u. 11}, 
            city={Budapest},
            postcode={1111}, 
            country={Hungary}}
\affiliation[angelina]{organization={Department of Mathematics, University of Michigan},%Department and Organization
            addressline={530 Church St}, 
            city={Ann Arbor},
            state={MI},
            postcode={48109}, 
            country={USA}}

%% Abstract
\begin{abstract}
%A dichotomy theorem on the parameterized complexity of the 3-uniform hypergraphicality problem is presented. 
We present a dichotomy theorem on the parameterized complexity of the 3-uniform hypergraphicality problem. Given $0<c_1\le c_2 < 1$, the parameterized 3-uniform Hypergraphic Degree Sequence problem, $\grparam$, considers degree sequences $D$ of length $n$ such that all degrees are between $c_1 {n-1 \choose 2}$ and $c_2 {n-1\choose 2}$ and it asks if there is a 3-uniform hypergraph with degree sequence $D$. We prove that for any $0<c_2< 1$, there exists a unique, polynomial-time computable $c_1^*$ with the following properties. For any $ c_1\in (c_1^*,c_2]$, $\grparam$ can be solved in linear time. In fact, for any $c_1\in (c_1^*,c_2]$ there exists an easy-to-compute $n_0$ such that any degree sequence $D$ of length $n\ge n_0$ and all degrees between $c_1 {n-1\choose 2}$ and $c_2 {n-1\choose 2}$ has a 3-uniform hypergraph realization if and only if the sum of the degrees can be divided by $3$. Further, $n_0$ grows polynomially with the inverse of $c_1-c_1^*$. On the other hand, we prove that for all $c_1<c_1^*$, $\grparam$ is NP-complete. Finally, we briefly consider an extension of the hypergraphicality problem to arbitrary $t$-uniformity. We show that the interval where degree sequences (satisfying divisibility conditions) always have $t$-uniform hypergraph realizations must become increasingly narrow, with interval width tending to $0$ as $t \rightarrow \infty$. 
% Finally, we briefly consider an extension of the hypergraphicality problem to arbitrary $t$-uniformity. We show that in higher uniformity, any interval of degree bounds within which a degree sequence has a realization (if it satisfies divisibility conditions) must become increasingly narrow, specifically with interval width tending to $0$ as $t \rightarrow \infty$.
% \textbf{TODO: add line about asymptotic result for $t$-uniformity.}
% (Arya) I made a command for the problem name: we can use \grprob and \grparam (for the version with c1, c2). See commands.tex.
\end{abstract}



%% Keywords

\end{frontmatter}

\section{Introduction}%% IN PROGRESS %%
The Degree Sequence Problem asks the following question: given a degree sequence, that is, a sequence of non-negative integers $D=(d_1,...,d_n),$ does there exist a simple graph $G=(V,E)$ such that the vertices have degrees corresponding to those in the sequence ($d(v_i) = d_i$ for $i=1,...,n$)? If the answer is yes, the graph $G$ is called a \textit{realization} of $D$, and $D$ is called a \textit{graphic} degree sequence. %As a start to answering this question, one can observe firstly that by the Handshake Lemma, which states that for any simple graph $G=(V,E), \sum_{i=1}^n d(v_i) = 2|E|,$ a necessary condition for $D$ having a realization is that $ \sum_{i=1}^n d_i$ is even. 
%Remark (Istvan): I would skip it. We might assume that the readers of the journal know the Handshake Lemma.

\vspace{0.5em}
This problem is one of the first solved problems in algorithmic graph theory. In 1955 and 1962, respectively, 
%\href{http://eudml.org/doc/19050}
Havel \cite{Havel} and 
%\href{https://www.jstor.org/stable/2098746}
Hakimi \cite{Hakimi} independently gave the same polynomial time algorithm to decide if a realization of a degree sequence $D$ exists. The algorithm is constructive; if a realization does exist, the algorithm gives such a realization. In 1960, 
%\href{https://www.renyi.hu/~p_erdos/1961-05.pdf}
Erd{\H o}s and Gallai \cite{ErdosGallai} also gave necessary and sufficient inequalities for the existence of a realization of a degree sequence. These inequalities can be easily checked in polynomial time; therefore, the graphicality problem for simple graphs is clearly in P.

\vspace{0.5em}
Hypergraphs are generalizations of simple graphs. A hyperedge $e \in E$ of a hypergraph $H=(V,E)$ is a non-empty subset of $V.$ A hypergraph is $k$-uniform if each edge is a subset of vertices of size $k.$ A hyperedge $e$ is incident with a vertex $v$ if $v \in e.$ The degree of a vertex is the number of hyperedges incident with it. With these definitions, the Degree Sequence Problem can be naturally generalized to hypergraphs: given a hypergraphic degree sequence $D=(d_1,...,d_n)$ and a positive integer $k,$ does there exist a $k$-uniform hypergraph $H=(V,E)$ such that the vertices have degrees corresponding to those in the sequence? We denote this problem by $\gr$ (the ``$k$-Uniform Hypergraphic Degree Sequence'' problem).
%Similarly to in the Degree Sequence Problem for simple graphs, one can observe that since the Generalized Handshake Lemma states that for any $k$-uniform hypergraph 
%$H=(V,E), \sum_{i=1}^n d(v_i) = k|E|,$ a necessary condition for $D$ to have a $k$-uniform realization is that $k | \sum_{i=1}^n d_i$.
%Remark (Istvan): I would skip it here, intriducing later
%
%Because of this, we will refer to degree sequences which have a realization as long as they satisfy the Generalized Handshake Lemma as "always graphic."
%Remark (Istvan): It is confusing. First, we should intriduce degree sequence classes. I'd like to skip it here, give a more detailed introduction of always graphicality later on.

\vspace{0.5em}
In 2018, Deza \emph{et al.} \cite{Dezaetal2018,Dezaetal2019} 
%\href{https://arxiv.org/pdf/1901.02272.pdf}{Deza et al.} 
proved that $\gr$  is NP-complete when $k=3.$ That is, it is NP-complete to decide if a 3-uniform hypergraph exists with a prescribed degree sequence. Given this result, it is natural to attempt to characterize the degree sequences for which \textsc{3Uni-HDS} can be solved in polynomial time. 

%The \href{https://arxiv.org/pdf/2312.00555}{2023 BSM Summer Research Paper}
\vspace{0.5em}
In 2023, Li and Mikl\'os \cite{limiklos} gave bounds $c_1$ and $c_2$ such that any degree sequence $D=(d_1,...,d_n)$ of length $n$ is graphic if $n$ is large enough, all degrees are between $c_1 n^2$ and $c_2 n^2$, and the sum of the degrees can be divided by $3$. The values are roughly $c_1= 0.03$ and $c_2=0.08$, corresponding to roughly $0.06 {n-1\choose 2}$ and $0.16 {n-1\choose 2}$.

\vspace{0.5em}
In this paper, we present a polynomial-time-computable threshold value $c_1^*$ for any $0<c_2<1$ with the following properties. (1) For sufficiently large $n$, a degree sequence of length $n$ is \emph{always graphic} when (a) the degrees are between $c_1{n-1\choose 2}$ and $c_2{n-1\choose 2}$ for some $c_1 > c_1^*$ and (b) the sum of the degrees is divisible by $3$; and (2) if the degrees are instead between $c_1 {n-1\choose 2}$ and $c_2{n-1\choose 2}$ for $c_1 < c_1^*$, the problem remains NP-complete. 
% (1) A degree sequence of length $n$ is always graphic when the degrees are between $c_1{n-1\choose 2}$ and $c_2{n-1\choose 2}$, $c_1 > c_1^*$, $n$ is large enough and the sum of the degrees can be divided by $3$; and (2) when there the degrees are between $c_1 {n-1\choose 2}$ and $c_2{n-1\choose 2}$, $c_1 < c_1^*$ the problem becomes NP-complete. 
% (Arya) Lightly edited the old version (commented out above) to make the sentence clearer and flow better, feel free to edit back if you disagree with the wording change
Our widest interval for always graphic degree sequences occurs when $c_1^{*}\approx 0.28$ and $c_2 \approx 0.72.$ Note that this is a substantially wider interval than the result in \cite{limiklos}.

\vspace{0.5em}
The first part of our result gives a lower bound on degrees in a degree sequence with a fixed length and maximum degree, % which guarantees that the sequence is always graphic (given that the sum of the degrees can be divided by $3$). 
which guarantees that the sequence is graphic, given that the sum of the degrees can be divided by $3$. Thus, the degree sequence class with the obtained lower and upper bounds is an always graphic degree sequence class. 
This lower bound falls out of a construction for the  realization, which is based on edge types and counting arguments. 

\vspace{0.5em}
The second part of our result claims that this lower bound is the lowest possible such bound. In fact, we prove that
for any $\varepsilon >0$, it is NP-complete to decide if a degree sequence of length $n$ with degrees between $c_1 {n-1\choose 2}$ and $c_2{n-1\choose 2}$, for $c_1 = c_1^*-\varepsilon$, has a $3$-uniform hypergraph realization.
%moving down by any $\epsilon > 0$ makes the problem NP-complete. We prove this result with an embedding argument. 
We use the result of Deza \emph{et al.} \cite{Dezaetal2019} 
%\href{https://arxiv.org/pdf/1901.02272.pdf}{Deza}'s result 
that the general problem is NP-complete for 3-uniform hypergraphs. In our proof, for any $\varepsilon >0$, we embed any arbitrary degree sequence $D_0$
%with degrees $\epsilon$ outside of the threshold values $c_1^{*}, c_2$ 
into a larger degree sequence $D$ within a very rigid construction so that $D$ is graphic if and only if $D_0$ is graphic. The length of $D$ is $n$, which is only polynomially larger than the length of $D_0$ for any fixed $\varepsilon >0$.
The embedding procedure depends on $\varepsilon$, and the degrees of $D$ are between  $c_1 {n-1\choose 2}$ and $c_2{n-1\choose 2}$, where $c_1 = c_1^*-\varepsilon$. 
%\href{https://arxiv.org/pdf/1901.02272.pdf}{Deza et al.} implies that the question of $D$ being graphic is NP-complete, so these two problems being equivalent means that the question of $D_0$ being graphic is also NP-complete, giving the desired result.

\vspace{0.5em}
Finally, we analyze how the question changes for the general problem of $t$-uniform hypergraphs. We claim that as $t \to \infty$, the lower bound threshold increases. That is, the always graphic bounds become increasingly narrower.
\section{Preliminaries}

We begin by introducing some basic definitions and notation related to hypergraphs. 

\begin{definition}{}
A hypergraph $H = (V, E)$ is a generalization of simple graphs. For all $e \in E,$ $e$ is a non-empty subset of $V.$ A hypergraph is $t$-uniform if for all $e \in E, e \in \binom{V}{t}.$
\end{definition}{}

\textit{Notation}: We will denote the empty (hyper)graph using the same notation for the empty set, i.e. $H = \emptyset$, where the meaning should be clear by context. When multiple (hyper)graphs are considered, we will use $V(H)$ to refer to the vertex set of a hypergraph $H$ for clarity. The induced subgraph (or subhypergraph) of $H$ on some subset of vertices $V' \subseteq V$ is denoted $H[V']$. The complete (hyper)graph on $n$ vertices is denoted by $K_n$. The bipartite (hyper)graph between two vertex parts $V_1, V_2 \subseteq V$ s.t. $V_1 \cap V_2 = \emptyset$ is denoted by $H[V_1, V_2]$, where in the hypergraph case this covers all edges in $H$ that are incident with at least one vertex in $V_1$ and at least one vertex in $V_2$.
%where all vertices belong to either $V_1$ or $V_2$. % In $3$-uniformity (or higher), $H[V_1, V_2, V_3]$ denotes the tripartite hypergraph between disjoint vertex sets $V_1, V_2, V_3$, i.e. covering all edges where all vertices are in either $V_1$, $V_2$, or $V_3$ and there is at least one vertex in each. 

\begin{definition}{}
\begin{enumerate}
    \item A hyperedge $e$ is incident with $v$ if $v \in e.$ The degree of a vertex $v$ of a hypergraph is the number of hyperedges incident with it, denoted by $d(v).$ The degree sequence of a hypergraph is the sequence of the degrees of its vertices, written as $(d_1,...,d_n),$ if $|V|=n$ and $d(v_i)=d_i\ \  \forall i=1,...,n.$
    \item A degree sequence $(d_1,\dots, d_n)$ is $k$-regular if $d_i=k$ for all $1\leq i\leq n$. A degree sequence is almost regular if for some $k$, $d_i=k$ or $d_i=k+1$ for all $1\leq i\leq n$. 
    % A degree sequence is $k$-regular if each degree is $k$. A degree sequence is almost regular if each degree is either $k$ or $k + 1$ for some $k.$
    \item Given a sequence $D$ of non-negative integers, we say that a hypergraph $H = (V, E)$ is a realization of $D$ if the sequence of the degrees of the vertices of $H$ is $D.$ If $D$ has a realization, then we say that $D$ is graphic.
\end{enumerate}
\end{definition}%{}
It is trivial to see that the following generalization of the Handshaking Lemma is true.
\begin{lemma}[Generalized Handshaking Lemma]\label{lem:handshaking}
Let $H = (V,E)$ be a $t$-uniform hypergraph. Then $\sum_{v\in V} d(v) \equiv 0$ (mod $t$).
\end{lemma}

We will consider a parametric decision problem on the graphicality of $3$-uniform hypergraphs. We start with a definition of degree sequence classes needed for parametrization. Notice that the Generalized Handshaking Lemma is considered in the definition.
\begin{definition}
Given $c_1, c_2 > 0,$ $\dcc$ denotes the class of $3$-uniform hypergraph degree sequences such that for each degree sequence $D \in \dcc$ of length $n$ the following holds:
\begin{enumerate}
    \item $\sum_{d \in D} d \equiv 0$  (mod $3$)
    \item $c_{1} \binom{n-1}{2} \leq d_i \leq c_{2}\binom{n-1}{2},$ for all $i=1,...,n$
\end{enumerate}
\end{definition}

Now, The parametric hypergraph degree sequence problem is the following.
\begin{definition}
    $\mathrm{\grparam}$:\\
    INPUT: Degree sequences $D = (d_1, \dots, d_n) \in \dcc$\\
    OUTPUT: ``Yes" if there exists a $3$-uniform hypergraph $H = (V, E)$ such that for all $i$ $d(v_i) = d_i$, and ``No” otherwise.
\end{definition}

We now introduce the operation known as a \textit{hinge flip} that will be a key tool in our analysis. Hinge flip operations were introduced first in approximating the permanent \cite{jerrumsinclair} and were recently popularized in network science \cite{rechneretal,amanatidiskleer,erdosetal}. We give the analogous operation for hypergraphs.

\begin{definition}[Hinge Flips]
% (cite https://arxiv.org/abs/2110.09068)
\begin{enumerate}
    \item A hinge flip operation on a realization $G = (V, E)$ of a degree sequence removes a(n) (hyper)edge $\{v_i\} \cup x \in E$ (for $x \in V$) and adds a(n) (hyper)edge $\{v_j\} \cup x \in E$, $v_i \neq v_j$. 
    \item The corresponding hinge flip operation on a degree sequence $D = (d_1,..., d_n)$ is an operation
    which decreases a $d_i$ in $D$ by 1 and increases a $d_j$ in $D$ by 1.
    \item If $d_i > d_j+1,$ we call it a balancing hinge flip; if $d_i = d_j+1$, we call it a neutral hinge flip; and otherwise, we call it a reverse hinge flip.
\end{enumerate}
\end{definition}

The following lemma and corresponding theorem were proved in \cite{limiklos}.

\begin{lemma}[\cite{limiklos}]
% (cite https://arxiv.org/pdf/2312.00555.pdf) 
Let $D$ be a graphic hypergraph degree sequence, and let $d_i, d_j \in D$ such that $d_i > d_j +1.$ Let $D'$ be the hypergraph degree sequence obtained from $D$ by subtracting 1 from $d_i$ and adding 1
to $d_j.$ Then any realization of $D$ has a balancing hinge flip operation yielding a realization of $D'$, and thus $D'$ is also a graphic hypergraph degree sequence.
\end{lemma}
\begin{theorem}[\cite{limiklos}]\label{thm:extremal-deg-seq}
% (also cite https://arxiv.org/pdf/2312.00555.pdf)
Let $D = (d_{\min},..., d_{\min}, d, d_{\max}, . . . , d_{\max})$ be a graphic hypergraph degree sequence on $n$ vertices with $d_{\min} \leq d \leq d_{\max}.$ Further let $D'$
be a hypergraph degree
sequence on $n$ vertices such that for all $d' \in D', d_{\min} \leq d' \leq d_{\max}$ and $\sum_{d \in D} d = \sum_{d' \in D'} d'.$ Then $D'$
is also graphic.
\end{theorem}

\begin{remk}
   From now on, when we consider the degree sequence  $D=(d_{\min},\ldots,d_{\min},d_{\min}+1,\ldots,d_{\min}+1,d_{\max}-1,\ldots,d_{\max}-1,d_{\max},\ldots,d_{\max})$ of length $n$ containing $k$ values equal to $d_{\max}-1$ or $d_{\max}$, and $n-k$ values equal to $d_{\min}$ or $d_{\min}+1$'s, we will refer to the $d_{\max}-1$ and $d_{\max}$'s as large degrees and to the $d_{\min}$ and $d_{\min}+1$'s as small degrees.
\end{remk}

% (Arya): Not sure how feasible or important it is, but we could also try to add 
%  (1) the nice parametrization (though that was for (f^*, d_max) pairs -- we should check: does it still transfer over to (c_1^*, c2) well), 
% and/or (2) plots showing how these values relate, for some concreteness

We now formally present the main result of this paper, stated in the following Dichotomy Theorem. Proving Part (I) is the focus of Section \ref{sec:always=graphic}, while proving Part (III) (from which Part (II) follows as a corollary) is the focus of Section \ref{sec:np-complete}.

\begin{theorem}[\textbf{Dichotomy Theorem}]\label{thm:dichotomy}
For any $0 < c_2 < 1$, there exists a \textit{unique}, \textit{polynomial-time computable} value $c_1^*$ ($0 < c_1^* < c_2$) such that the following holds: 
\begin{enumerate}
\item [(I)] $\forall c_1>c_1^{*}$, $\grparam$ can be solved in linear time. In fact, $\exists \ n_0 = O(\poly(\frac{1}{c_1-c_1^{*}}))$ such that $\forall n \geq n_0$, any $n$-length degree sequence $D = (d_1, \dots, d_n) \in \dcc$ has a $3$-uniform hypergraph realization.
%of length $n \geq n_0$ with degrees in these bounds, $D$ 

%if and only if $3 | \sum_{i=1}^n d_i$.
% \vspace{0.5em} Furthermore, the proof is \textit{constructive}, and a realization can be obtained in $\poly(n)$ time (when one exists).}
\item [(II)] $\forall c_1 < c_{1}^{*}$, $\grparam$ is NP-complete. 
\item [(III)] In fact, $\forall \varepsilon > 0$ %$\exists \lambda$, such that 
the decision problem over the class of degree sequences $D = (d_1, \dots, d_n)$ of length $n$ where $\forall d_i \in D$, $c_1^* \binom{n-1}{2} - %\lambda
n^{1+\varepsilon} \leq d_i \leq c_2 \binom{n-1}{2}$, is NP-complete. 
\end{enumerate}
\end{theorem}

%\section{Linearly bounded always graphic sequences in $3$-uniformity}
\section{Linearly bounded always graphic 3-uniform hypergraph degree sequences}
\label{sec:always=graphic}

In this section, we begin by defining a class of degree sequences called \emph{critical degree sequences}. Then we will show in Lemma \ref{lem:critical-degree} that each degree sequence in this class has a $3$-uniform hypergraph realization, through a construction that motivates the specific numerical expressions in Definition~\ref{def:critical-degree-class}. Critical degree sequences plays a central role in proving the dichotomy theorem. We transform critical degree sequences into degree sequences $\dmin, \ldots, \dmin, d_{int}, \dmax,\ldots,\dmax$ with certain $\dmin$, $\dmax$ values and $\dmin \le d_{int}\le \dmax$, and we show that these transformed degree sequences are also graphic (Lemma~\ref{lem:g-int-graphic}). This will be the key to obtain always graphic degree sequence classes via Theorem~\ref{thm:extremal-deg-seq}. Further, in Section~\ref{sec:np-complete}, we will use critical degree sequences in an embedding process to prove the NP-completeness part of our dichotomy theorem.
%the following definition. 

\begin{definition}\label{def:critical-degree-class}
    The \emph{critical degree sequence class}, $\mathcal{D}_{crit}$, contains degree sequences that are each parameterized as follows by parameters $n,k,\dmax$.
    %is a parameterized degree sequence class with parameters $n,k,\dmax$. 
    The degree sequence $D(n,k,\dmax)\in \mathcal{D}_{crit}$ has $n$ degrees. The parameters must satisfy $\dmax \le {n-1\choose 2}$ and $k\in\{1,2,\ldots n\}$, with the following additional restrictions:
    \begin{itemize}
    \item If $k\dmax \equiv 1\pmod{3}$ then $k \le n-2$.
    \item If $k\dmax \equiv 2\pmod{3}$ then $k \le n-1$.
    \item If $\binom{k-1}{2}<\dmax\leq \binom{k-1}{2}+(n-k)(k-1)$ and $k\left(\dmax - \binom{k-1}{2}\right)\equiv 1\pmod{2}$ then $k \le n-2$.
    \end{itemize}
    Then $D(n,k,\dmax)$ contains $k$ $\dmax$ degrees as large degrees and $n-k$ small degrees. The small degrees are the following:
    \begin{itemize}
        \item If $\dmax \le {k-1\choose 2}$ then there are $2k\dmax \pmod{3}$ degrees of $1$ and $n-k-\left(2k\dmax \pmod{3}\right)$ degrees of $0$.
        \item If ${k-1\choose 2} < \dmax \le {k-1\choose 2} +(n-k)(k-1)$ then let $s := \left\lfloor\frac{k\left(\dmax-{k-1\choose 2}\right)}{2}\right\rfloor + 2\times\left(k\left(\dmax-{k-1\choose 2}\right) \pmod{2}\right)$. There are $s \pmod{n-k}$ degrees of $\ceil{\frac{s}{n-k}}$  and $n-k-\left(s\pmod{n-k}\right)$ degrees of $\floor{\frac{s}{n-k}}$.
        \item If ${k-1\choose 2} +(n-k)(k-1) < \dmax$ then let $s:= {k\choose 2}(n-k) + 2k\left(\dmax - {k-1\choose 2} - (n-k)(k-1)\right)$. There are $s \pmod{n-k}$ degrees of $\ceil{\frac{s}{n-k}}$ degrees and $n-k-\left(s\pmod{n-k}\right)$ degrees of $\floor{\frac{s}{n-k}}$.
    \end{itemize}
\end{definition}

Before we prove that each degree sequence in $\mathcal{D}_{crit}$ has a $3$-uniform hypergraph realization, we give the following definition to classify edge types, which will be useful for many arguments throughout the remainder of the paper. 

\begin{definition}[3L, 2L1S, 1L2S, 3S edge types] \label{def:hyperegde-types}
Consider a $3$-uniform hypergraph $H=(V_L\sqcup V_S,E)$, where $\sqcup$ denotes the disjoint union of vertex sets. Call $V_L$ the large degree vertices, and $V_S$ the small degree vertices. Then, we can define the following edge types:
\begin{itemize}
    \item An edge $e\in E$ is a 3L edge if $e\subseteq V_L$. 
    \item An edge $e\in E$ is a 2L1S edge if $|e\cap V_L|=2$ and $|e\cap V_S|=1$. 
    \item An edge $e\in E$ is a 1L2S edge if $|e\cap V_L|=1$ and $|e\cap V_S|=2$. 
    \item An edge $e\in E$ is a 3S edge if $e\subseteq V_S$. 
\end{itemize}   
\end{definition}

%\noindent 

\begin{lemma}\label{lem:critical-degree}
    Each degree sequence $D(n,k,\dmax) \in \mathcal{D}_{crit}$ has a $3$-uniform hypegraph realization.
\end{lemma}
\begin{proof}
    
We can construct a hypergraph realization $H= (V_L\sqcup V_S, E)$ of $D(n,k,\dmax)$ in the following way:\\
\\
First, arbitrarily add edges of specified types according to the three cases below: 
\begin{itemize}
    \item If $\dmax \leq \binom{k-1}{2}$, add $\floor{\frac{k\dmax}{3}} $ 3L edges. Note that indeed $\floor{\frac{k\dmax}{3}} \leq \binom{k}{3}$, the total number of possible 3L edges, since $\dmax \leq \binom{k-1}{2}$. If $k\dmax \equiv 1\pmod{3}$ or $k\dmax \equiv 2\pmod{3}$, add one 1L2S or one 2L1S edge, respectively. Further, $k\dmax$ can be congruent with $2$ modulo $3$ only if $k\le n-1$ and can be congruent with with $1$ modulo $3$ only if $k\le n-2$. Therefore these 1L2S or 2L1S edges are available, that is, there are sufficient small degree vertices.
    
    \item If $\binom{k-1}{2}<\dmax\leq \binom{k-1}{2}+(n-k)(k-1)$, add all $\binom{k}{3}$ 3L edges and $\floor{ \frac{k\left(\dmax - \binom{k-1}{2}\right)}{2}}$ 2L1S edges. Note that indeed $\floor{\frac{k\left(\dmax - \binom{k-1}{2}\right)}{2}} \leq \binom{k}{2}(n-k),$ the total number of possible 2L1S edges, since $\dmax \leq \binom{k-1}{2}+(n-k)(k-1).$ If $k\left(\dmax - \binom{k-1}{2}\right)\equiv 1\pmod{2}$ then
    %this means that there is one large degree vertex left which has degree $\dmax-1:$ 
    add one 1L2S edge. We have not added any 1L2S edges yet, and due to the restriction on $k$, this edge should be available, that is, there are sufficient small degree vertices.

    \item If $\binom{k-1}{2}+(n-k)(k-1)<\dmax$, add all $\binom{k}{3}$ 3L edges, all $\binom{k}{2}(n-k)$ 2L1S edges, and $k\left(\dmax-\binom{k-1}{2}-(n-k)(k-1)\right)$ 1L2S edges. Since $\dmax \leq \binom{n-1}{2}$ and one can prove that $\binom{n-1}{2}-\binom{k-1}{2}-(n-k)(k-1) = \binom{n-k}{2},$ we have that $k\left(\dmax-\binom{k-1}{2}-(n-k)(k-1)\right) \leq k\binom{n-k}{2},$ the total number of possible 1L2S edges.
\end{itemize}
Then, perform balancing hinge flips as follows:
\begin{enumerate}
    \item Let $v_1$ be a vertex with degree $\max_{v \in V_L} d(v)$ and let $v_2$ be a vertex with $\min_{v \in V_L} d(v)$. % Take $v_1\in V_L$ such that $d(v_1)=\max \{d(v):v\in V_L\}$ and $v_2\in V_L$ such that $d(v_2)=\min \{d(v):v\in V_L\}$. 
    If $d(v_1) > d(v_2)+1$, perform a balancing hinge flip. Repeat this step while there exists degrees $d(v_1) > d(v_2)+1$. We claim that this procedure arrives to a regular degree sequence on $V_L$. Indeed, note that $\sum_{v\in V_L} d(v) = |V_L| \dmax = k \dmax$. Therefore, if there exists a $d(v_1) > \dmax$ then there also exists a $d(v_2) < \dmax$. Also, if there exists a $d(v_1) < \dmax$ then there exists a $d(v_2) > \dmax$. Performing a balancing hinge flip on these degrees decreases $\sum_{v\in V_L} |d(v)-\dmax|$. Therefore, in finite number of steps, we will arrive to $\max_{v \in V_L} d(v) = \min_{v \in V_L} d(v) = \dmax$, that is, the degree sequence segment of $V_L$ is regular.
    \item Let $u_1$ be a vertex with degree $\max_{v \in V_S} d(v)$ and let $u_2$ be a vertex with degree $\min_{v \in V_S} d(v)$.  
    % Take $u_1\in V_S$ such that $d(u_1)=\max \{d(u):u\in V_S\}$ and $u_2\in V_S$ such that $d(u_2)=\min \{d(u):u\in V_S\}$. 
    If $d(u_1)-d(u_2)\geq 2$, perform a balancing hinge flip. Repeat this step while $d(u_1)-d(u_2)\geq 2$. It is easy to see that this procedure arrives to an almost regular degree sequence on $V_S$.
\end{enumerate}
To finish the proof, we are going to show that the almost regular degrees on the small degree vertices $V_S$ in $H$ are indeed the ones given in the definition of $D(n,k,\dmax)$ (Definition \ref{def:critical-degree-class}).
\begin{itemize}
    \item If $\dmax \le {k-1 \choose 2}$ then there are at most $2$ degrees of $1$ in $V_S$, all other degrees are $0$. If $k\dmax \equiv 0 \pmod{3}$, then all small degrees are $0$. If $k\dmax \equiv 1 \pmod{3}$, then there are $2$ degrees of $1$ since one 1L2S hyperedge is added. Indeed $2\times 1 \equiv 2 \pmod{3}$. If $k\dmax \equiv 2 \pmod{3}$, then there are $1$ degree of $1$ since one 2L1S hyperedge is added. Indeed, $2\times 2\equiv 1 \pmod{3}$.
    \item If ${k-1\choose 2} < \dmax \le {k-1\choose 2} + (n-k)(k-1)$, then the sum of degrees of all vertices in $V_S$ is $s:=\left\lfloor\frac{k\left(\dmax-{k-1\choose 2}\right)}{2}\right\rfloor + 2\times\left(k\left(\dmax-{k-1\choose 2}\right) \pmod{2}\right)$. Indeed, $\left\lfloor\frac{k\left(\dmax-{k-1\choose 2}\right)}{2}\right\rfloor$ 2L1S hyperedges are added and one 1L2S hyperedge is added if $k\left(\dmax-{k-1\choose 2}\right) \equiv 1 \pmod{2}$. This sum has to be distributed almost regularly among $n-k$ degrees. Then indeed there are $s \pmod{n-k}$ degrees of $\left\lceil\frac{s}{n-k}\right\rceil$ degrees and $n-k-\left(s\pmod{n-k}\right)$ degrees of $\left\lfloor\frac{s}{n-k}\right\rfloor$.
    \item If ${k-1\choose 2} +(n-k)(k-1) < \dmax$ then the sum of the small degree vertices is $s:= {k\choose 2}(n-k) + 2k\left(\dmax - {k-1\choose 2} - (n-k)(k-1)\right)$. Indeed, all 2L1S hyperedges are added and $k\left(\dmax - {k-1\choose 2} - (n-k)(k-1)\right)$ 1L2S hyperedges. This sum has to be distributed almost regularly among $n-k$ degrees. Then indeed there are $s \pmod{n-k}$ degrees of $\left\lceil\frac{s}{n-k}\right\rceil$ degrees and $n-k-\left(s\pmod{n-k}\right)$ degrees of $\left\lfloor\frac{s}{n-k}\right\rfloor$.
\end{itemize} \end{proof}

We will refer to a(n arbitrary) realization of a critical degree sequence, guaranteed to exist by Lemma \ref{lem:critical-degree}, as a \textit{critical hypergraph}. 

\vspace{0.5em}
We now proceed to define the following functions $f(n, k, \dmax)$ (Definition \ref{def:f}) and $g(n, k, \dmax)$ (Definition \ref{def:g}) of the parameters $n, k, \dmax$ based on critical degree sequences and critical hypergraphs. These functions define lower bounds on the always graphic interval given $k$ and $\dmax$, and their limits as $n$ becomes large will be the key ingredients that define the critical threshold value $c_1^*$.     

\begin{definition}\label{def:f}
Define $f_0(n,k,\dmax)$ as the average degree of the small degrees in the critical degree sequence $D(n,k,\dmax)$. Then define $f(n,k,\dmax) := \ceil{f_0(n, k, \dmax)}$ and $\fopt := \max_k \f$. Further, let $k^*(n,\dmax) := \arg\max_k \f$
\end{definition}

\begin{lemma}\label{lem:theta-k}
    Let $c_2 \in (0,1)$ be an arbitrary real number. Then there exists an $\varepsilon > 0$ and $n_0\in \mathbb{N}$ such that for all $n>n_0$, $\frac{k^*(n,c_2{n-1\choose2})}{n} >\varepsilon$ and $\frac{k^*(n,c_2{n-1\choose 2})}{n} < 1-\varepsilon$. Further, for any $n>n_0$,${k^*(n,c_1{n-1\choose 2})-1 \choose 2} < c_2{n-1\choose 2}$.
\end{lemma}
\begin{proof}
    For sufficiently large $n$, there exists a $k$ such that ${k-1\choose2}<c_2{n-1\choose 2} \le {k-1\choose2}+(n-k)(k-1)$. For any such $k$, $f_0(n,k,c_2{n-1\choose2}) = \Omega(n^2)$, further, both $k$ and $n-k$ are $\Omega(n)$. On the other hand, whenever $k = o(n)$ or $n-k = o(n)$, $f_0(n,k,c_2{n-1\choose 2}) = o(n^2)$. 

    To prove the second statement of the lemma, simply observe that for all $k$ such that ${k-1 \choose 2} \ge c_2{n-1\choose 2}$, $f_0(n,k,c_2{n-1\choose2}) = o(n)$.
\end{proof}

\begin{definition}\label{def:g}
Fix $n, k, \dmax$ as before. Define $g(n, k, \dmax) := \f + \ceil{\frac{2(\dmax - \f)}{n-k-1}}$. Define $\gopt := \max_k g(n, k, \dmax)$.
\end{definition}

\begin{lemma}\label{lem:g-int-graphic}
Fix $n,k,\dmax$. Define the degree sequence $D=(\dm,\dots,\dm,\dint,\dmax,\dots,\dmax)$, where $k$ degrees are $\dmax$, $n-k-1$ degrees are $\dmin$, and $\dmin\leq \dint \leq \dmax \leq {n-1\choose 2}$. If $\sum_{d\in D}d\equiv 0\pmod{3}$ and $\dm \in [g(n,k,\dmax),\dmax]$, then $D$ is graphic. 
\end{lemma}
\begin{proof}
Consider a critical hypergraph with parameters $n$, $k$, and $\dmax$,
%critical degree sequence $D(n,k,\dmax)$ and a realization of it, 
$H_0=(V_L\sqcup V_S, E)$. Fix $v_{int}$, where $v_{int}\in V_S$ and $d(v_{int})=\max \{d(v):v\in V_S\}$. While $d(v_{int})< \dint$, add an edge $e$ not present in the current realization such that $v_{int}\in e$. This is possible since $d(v_{int})<\dint \le\binom{n-1}{2}$. When this process terminates, call the resulting graph $H'=(V'_L\sqcup V'_S,E')$. Let $s':= \sum_{v\in V'_L\sqcup V'_S} d(v)$, and let $s := \sum_{d\in D}d=(n-k-1)\dm + \dint + k \dmax$. 

\vspace{0.5em}
Observe that $s-s'\equiv 0\pmod{3}$. Furthermore, we claim that it is non-negative. We have that $s' \leq (n-k-1) \cdot f(n, k, \dmax) + k \dmax + \dint + 2 \cdot m$ where $m$ is the number of edges added to create $H'$ from $H_0$. Observe that $m = \dint - f(n, k, \dmax) \leq \dmax - f(n, k, \dmax)$, since in $H_0$, $d(v_{int}) = f(n, k, \dmax)$. Thus $s - s' \geq (n-k-1) \cdot (\dmin - f(n,k,\dmax)) - 2 \cdot (\dmax - f(n, k, \dmax))$ where $\dmin \geq g(n,k,\dmax)$, and then by definition of $g(n, k, \dmax)$ it follows that $s - s' \geq 0$. 

\vspace{0.5em}
Thus we can add $(s-s')/3$ arbitrary hyperedges to $H'$. We keep calling this hypergraph $H' = (V'_L \sqcup V'_S, E')$. 
%\textbf{(TODO: fix, can't quite do this ``arbitrarily'')}. 
Then do the following balancing hinge-flips:
\begin{enumerate}
    \item While there is a $v\in V'_L$ such that $d(v) >\dmax$, let $u\in V'_S\setminus\{v_{int}\}$ be a vertex with minimal degree. We claim that $d(v) > d(u)+1$.
    Indeed, each vertex in $V'_L$ has a degree at least $\dmax$, $d(v) > \dmax$ and $v_{int}$ has a degree at least $\dint$. If the smallest degree in $V'_S$ were at least $\dmax$, then it would contradict that the sum of the degrees is $(n-k-1)\dmin + \dint + k \dmax$. Therefore, $d(u) < \dmax$.
    Do a balancing hinge-flip between $u$ and $v$.
    \item While $d(v_{int}) > \dint$ (this can be happen due to adding $(s-s')/3$ hyperedges to $H'$), let $u\in V'_S\setminus\{v_{int}\}$ be a vertex with minimal degree. Similarly to the previous point, it is easy to see that $d(u) +1 < d(v_{int})$. Do a balancing hinge-flip between $u$ and $v_{int}$.  
    \item While there are two vertices $u,v \in V_S \setminus\{v_{int}\}$ with $d(u)-d(v) \ge 2$, do a balancing hinge-flip between $u$ and $v$. Since the average degree on $V'_S \setminus \{v_{int}\}$ is $\dmin$, this procedure terminates in a regular degree sequence on $V'_S \setminus \{v_{int}\}$.
\end{enumerate}
The resulting hypergraph is a realization of $D$. \end{proof} % (Arya) could add a quick line of justification why, if needed, later -- V_L and v_int are correct now, and then V_S is almost regular with the correct sum => in fact regular dmin

Putting Lemma \ref{lem:g-int-graphic} together with Theorem \ref{thm:extremal-deg-seq} yields the following lemma on the graphicality of degree sequences with degrees between $g^*(n, \dmax)$ and $\dmax$. 

\begin{lemma}\label{lem:always-graphic-fix-n}
    Let $D$ be a degree sequence on $n$ vertices. Let $\dmax$ be its largest degree and $\dmin$ be its smallest degree. If the sum of the degrees in $D$ can be divided by $3$ and $\dmin$ is at least $g^*(n,\dmax)$, then $D$ is graphic.
\end{lemma}
\begin{proof}
    Given $\dmin$ and $\dmax$, observe that we can find $k$ and $d_{int}$ such that
    $$
    (n-k-1)\dmin + d_{int} + k \dmax = \sum_{d_i\in D} d_i,
    $$
    and $\dmin\le d_{int} < \dmax$.
    We know from Definition~\ref{def:g} and Lemma~\ref{lem:g-int-graphic} that the degree sequence containing $n-k-1$ degrees of $g^*(n,\dmax)$, one degree $d_{int}$, and $k$ degrees of $\dmax$ is graphic. Then, by Theorem \ref{thm:extremal-deg-seq}, $D$ is also graphic. 
    % $D'$ degree sequence containing $n-k-1$ degrees of $g^*(n,\dmax)$ a degree $d_{int}$ and $k$ degrees of $\dmax$ is graphic. 
\end{proof}

% \subsection{$c_1^*$ proofs}

% {\color{black!50} \textbf{8/31} (Arya): I think I have a proof idea for all the main $c_1^*$ related claims. The main challenge with \textit{}the rounding/continuous relaxation approach that Prof. Miklos proposed seemed to be that both the piecewise function values and the piecewise case bounds depend on $n$ and $k$ (which are being modified/discretized). My proof is based on a ``two-step rounding'' approach (first round the values, then round the bounds, and prove that each step changes the pointwise value by $O(n)$). Then I think there needs to be another rounding step for turning $k$ and $\dmax$ into continuous values, still needs to be fleshed out. One important thing that this proof relies on is that $n-k^* = \Theta(n)$, which I believe we have proved but should ensure it all fits together. A lot of this writing is currently work in progress, notation might not align well / claims might be rough around the edges -- please check.}

With Lemma \ref{lem:always-graphic-fix-n}, we are now close to obtaining Part (I) of our dichotomy theorem. What remains is to use the $g^*(n, \dmax)$ function, which depends on $n$, to derive a critical value $c_1^*$ that does not depend on $n$ (which we will show has the desired properties). This is accomplished in the following definition and subsequent technical lemmas \ref{lem:discrete-continuous-f0-g} and \ref{lem:conv-f-star-g-star}. 

\begin{definition}[Critical value $c_1^*$]\label{def:c1star}
%% add bounds on c_2 explicitly?
Fix $c_2\in (0,1)$. The critical value $c_1^*$ corresponding to $c_2$ is defined as
$$c_1^*(c_2) = \max_{\alpha \in (0, 1)} C\left(\alpha, \frac{c_2}{2}\right)$$
%% add bounds on d explicitly?
where $C(\alpha, d)$ is given by 

\begin{equation*}
    C(\alpha, d) = \begin{cases} 
      0 & d \leq \frac{\alpha^2}{2} \\
    %  \frac{\alpha}{1-\alpha} \left(\frac{2d-\alpha^2}{4}\right) & 
       \frac{\alpha}{1-\alpha} \left(\frac{2d-\alpha^2}{2}\right) & 
      \frac{\alpha^2}{2} < d \leq \alpha(1-\frac{\alpha}{2}) 
      % \alpha^2 < d \leq \alpha(2-\alpha) 
      \\ 
    %  \frac{\alpha}{1-\alpha}(2d-\alpha^2) - \frac{3\alpha^2}{2} & d > \alpha(1-\frac{\alpha}{2})
      \frac{2\alpha}{1-\alpha}(2d-\alpha^2) - 3\alpha^2 
      & d > \alpha(1-\frac{\alpha}{2})
      % & d >  \alpha(2-\alpha)
   \end{cases}
\end{equation*}

We also define $$\alpha^* := \arg\max_{\alpha \in (0, 1)} C\left(\alpha, \frac{c_2}{2}\right).$$
\end{definition}

\begin{observation}\label{obs:never-zero}
For any $0<c_2<1$, $c_1^*(c_2) > 0$ and $0<\alpha^* <1$.
\end{observation}
\begin{proof}
    Consider any $\alpha \in (1-\sqrt{1-c_2},\sqrt{c_2})$ (it is easy to see that this interval is not empty). Then $C\left(\alpha,\frac{c_2}{2}\right) > 0$ and thus $c_1^*(c_2)>0$. Further, for all $c_2 \in (0,1)$, $C\left(0,\frac{c_2}{2}\right)=0$ and $C\left(1,\frac{c_2}{2}\right) = 0$, thus $\alpha^* \in (0,1)$.
\end{proof}

\begin{lemma}\label{lem:discrete-continuous-f0-g}
    For any $\alpha \in (0,1)$ and $c_2 \in (0,1)$,
    $$
    \lim_{n\rightarrow \infty} \frac{f_0\left(n,\floor{\alpha n},\floor{c_2 {n-1\choose 2}}\right)}{{n-1\choose 2}} 
    =  \lim_{n\rightarrow \infty} \frac{g\left(n,\floor{\alpha n},\floor{c_2 {n-1\choose 2}}\right)}{{n-1\choose 2}} = C\left(\alpha, \frac{c_2}{2}\right).
    $$
    Further, the convergence is polynomially fast. That is, there exists a universal polynomial $poly$ such that $\forall \alpha \in (0,1), c_2 \in (0,1)$  and $\forall\varepsilon>0$,  there exists an $n_0$ such that $n_0 = O(\poly(\frac{1}{\varepsilon}))$ and for all $n\ge n_0$,
    $$
    \left|\frac{f_0\left(n,\floor{\alpha n},\floor{c_2 {n-1\choose 2}}\right)}{{n-1\choose 2}} - C\left(\alpha, \frac{c_2}{2}\right) \right| \le \varepsilon
    $$
    and
    $$
    \left|\frac{g\left(n,\floor{\alpha n},\floor{c_2 {n-1\choose 2}}\right)}{{n-1\choose 2}} - C\left(\alpha, \frac{c_2}{2}\right) \right| \le \varepsilon
    $$
\end{lemma}

\noindent The proof -- although straightforward -- is quite technical, and therefore, it is given in \ref{app:3}.

\begin{lemma}\label{lem:conv-f-star-g-star}
    For any $0<c_2<1$,
    $$
    \lim_{n\rightarrow \infty}  \frac{f^*(n,\floor{c_2 {n-1\choose 2}})}{{n-1\choose 2}} =
    \lim_{n\rightarrow \infty}  \frac{g^*(n,\floor{c_2 {n-1\choose 2}})}{{n-1\choose 2}} =
    c_1^*(c_2).
    $$
    Further, the convergence is polynomially fast. That is, $\forall c_2 \in (0,1)$ and $\forall \varepsilon > 0$, $\exists n_0$ such that $n_0 = O(\poly(\frac{1}{\varepsilon}))$ and for all $n \ge n_0$
    $$
    \left|\frac{f^*(n,\floor{c_2 {n-1\choose 2}})}{{n-1\choose 2}}- c_1^*(c_2) \right| \le \varepsilon
    $$
    and
    $$
    \left|\frac{g^*(n,\floor{c_2 {n-1\choose 2}})}{{n-1\choose 2}}- c_1^*(c_2) \right| \le \varepsilon
    $$
\end{lemma}

\noindent The proof is also quite technical, and therefore it is given in \ref{app:3}.

  % The speed of convergence follows from $\sup_{\alpha\in (0,1)}\left\{n(\alpha)\right\} = O(poly(\frac{1}{\varepsilon}))$ due to Lemma~\ref{lem:discrete-continuous-f0-g}.

\vspace{0.5em}
\noindent For the hardness part of the main theorem, we need the following lemma.
 \begin{lemma}\label{lem:o-n-to-1-plus-epsilon}
 For any $0<c_2<1$, $\lambda >0$ 
 and  $\varepsilon > 0$, there exists an $n_0$ such that for all $n\ge n_0$,
 $$
 c_1^*(c_2){n-1\choose 2} - \lambda 
 n^{1+\varepsilon} \le f^*(n,\left\lfloor c_2{n-1\choose 2}\right\rfloor).
 $$    
 \end{lemma}

 \noindent The technical proof is again given in \ref{app:3}.

\vspace{0.5em}
We are now ready to prove Part (I) of the Dichotomy Theorem, which follows in a straightforward manner from the pieces built up so far.

% somehow indicate / note that this is (or equivalently proves) Dichotomy Theorem (I)
\begin{theorem}{\normalfont(Part (I) of \hyperref[thm:dichotomy]{Dichotomy Theorem})}
Let $0<c_1\le c_2 <1$ be real values such that $c_1 > c_1^*(c_2)$. Then there exists an $n_0 = O(\poly(\frac{1}{c_1^*(c_2)-c_1}))$ such that for any $n \ge n_0$, any degree sequence $D$ of length $n$ with degrees between $c_1 {n-1 \choose 2}$ and $c_2 {n-1\choose 2}$ has a $3$-uniform hypergraph realization if and only if the sum of its degree can be divided by $3$. In particular, any degree sequence $D$ of length $n \geq n_0$ in $\dcc$ is graphic.
\end{theorem}
\begin{proof}
    Let $\varepsilon := c_1-c_1^*(c_2)$.
    By Lemma~\ref{lem:conv-f-star-g-star}, there exists an $n_0 = O(\poly(\frac{1}{\varepsilon})) = O(\poly(\frac{1}{c_1 - c_1^*(c_2)}))$ such that for all $n\ge n_0$,
    $$
    \frac{g^*(n,\floor{c_2 {n-1\choose 2}})}{{n-1\choose 2}} - c_1^*(c_2) \le \varepsilon,
    $$
    that is,
    $$
    g^*\left(n,\floor{c_2 {n-1\choose 2}}\right) \le (c_1^*(c_2)+\varepsilon) {n-1\choose 2} = c_1{n-1\choose 2}.
    $$ 
    %Further, since $\frac{g^*(n,\floor{c_2{n-1\choose 2}}}{{n-1\choose 2}}$ can be lower and upper bounded by rational polynomial functions both converging to $c_1^*(c_2)$, $n_0 = O(poly(\frac{1}{c_1^*(c_2)-c_1}))$. 
    Thus, it follows from Lemma~\ref{lem:always-graphic-fix-n} that for any $n\ge n_0$, any degree sequence $D$ of length $n$ with degrees between $c_1 {n-1 \choose 2}$ and $c_2 {n-1\choose 2}$ has a $3$-uniform hypergraph realization if and only if the sum of its degrees can be divided by $3$.
\end{proof}

\begin{figure}
    \centering
    \includegraphics[width=0.7\linewidth]{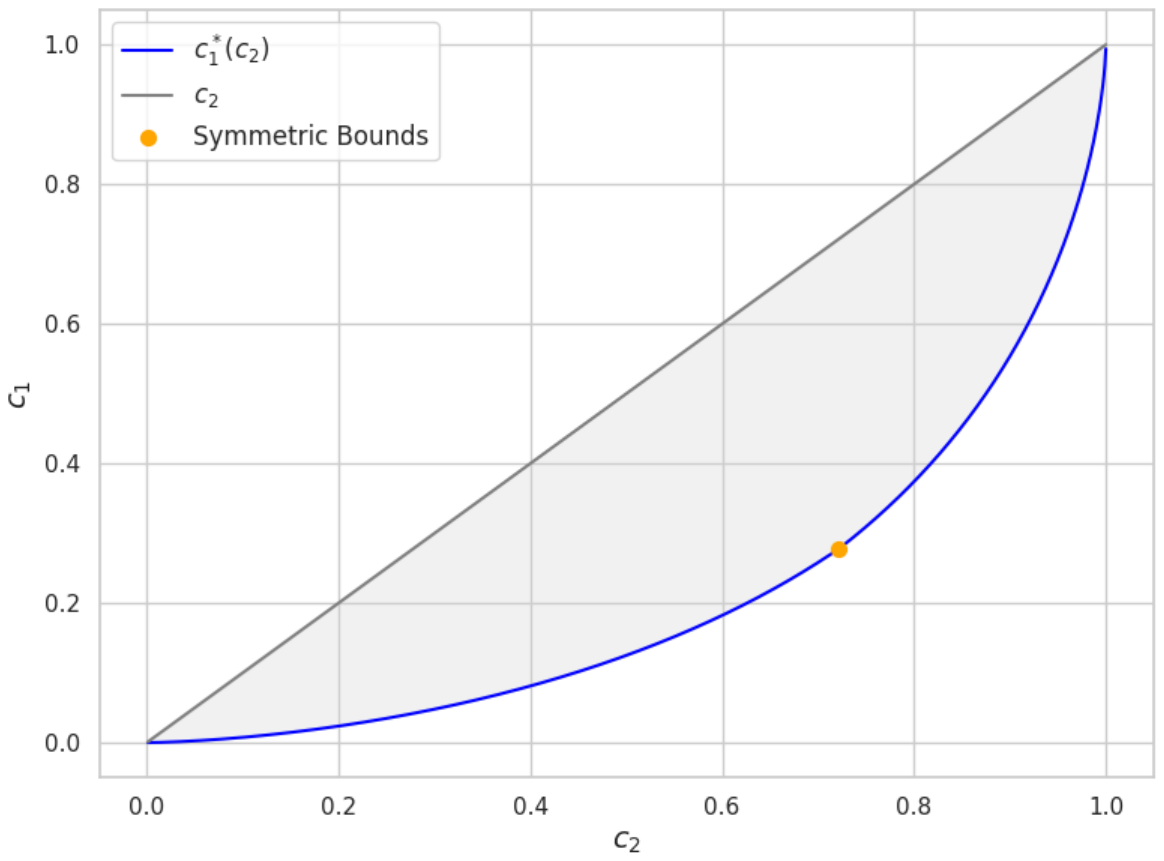}
    \caption{Plot of $c_1^*(c_2)$ for all $c_2 \in (0, 1)$. The shaded region indicates all $c_1$ s.t. $c_1 \geq c_1^*(c_2)$, i.e. the region where $\grparam$ is easily solvable according to Part (I) of the Dichotomy Theorem (Theorem \ref{thm:dichotomy}). The orange point indicates the symmetric bounds case ($c_2 \approx 0.721934$, $c_1^* \approx 0.278066$), which are empirically the widest bounds obtained.}
    \label{fig:c1star-c2-plot}
\end{figure}

One might ask which $(c_1^*(c_2),c_2)$ interval is the widest, that is, when $c_2-c_1^*(c_2)$ is maximal. Empirical results suggest that it happens when $c_1^*(c_2) = 1-c_2$, and  $\frac{\alpha^2}{2} < \frac{c_2}{2} \le \alpha\left(1-\frac{\alpha}{2}\right)$ (see also Figure~\ref{fig:c1star-c2-plot}).
%Therefore, by the definition of $c_1^*$ (Definition \ref{def:c1star}), it 
This symmetric case can be expressed as the unique solution between $0$ and $1$ for $c_2$ of the following equation system:
$$
\frac{\alpha}{1-\alpha} \frac{c_2-\alpha^2}{2} = 1- c_2,
$$
$$
\frac{d}{d\alpha}\left(\frac{\alpha}{1-\alpha} \frac{c_2-\alpha^2}{2}\right) = \frac{1}{2}\frac{(c_2-3\alpha^2)(1-\alpha)+\alpha(c_2-\alpha^2)}{(1-\alpha)^2} = 0.
$$
The approximate value for $c_2$ is $0.721934$, the corresponding $c_1^*(c_2)$ is approximately $0.278066$. The corresponding $\alpha$ value, which represents the so-called ``critical density'' of large degree vertices (i.e. the fraction of vertices in the critical degree sequence which are large degree), is approximately $0.652704$. It is worth mentioning that the corresponding widest interval for simple graphs is the $\left(\frac{1}{4},\frac{3}{4}\right)$ interval with $0.5$ being the critical density \cite{ems2024}.

%{\bf TODO: explain here that $c_1^*(c_2)$ is never given on the Case 0 part.} \textcolor{red}{covered by observation 3.9 below?}

Recall from Observation~\ref{obs:never-zero} that $c_1^*(c_2)$ is taken as the maximum of $\frac{\alpha}{1-\alpha} \frac{c_2-\alpha^2}{2}$ or $\frac{2\alpha}{1-\alpha}(c_2\alpha^2) - 3\alpha^2$, or explicitly at $c_2 = \alpha(2-\alpha)$ where $C\left(\alpha,\frac{c_2}{2}\right)$ is not differentiable. When $c_2 = \alpha(2-\alpha)$, then $C\left(\alpha,\frac{c_2}{2}\right) = (1-\sqrt{1-c_2})^2$.
% Recall Observation~\ref{obs:never-zero} that $c_1^*(c_2)$ is taken as the maximum of $\frac{\alpha}{1-\alpha} \frac{c_2-\alpha^2}{2}$ or the maximum of $\frac{2\alpha}{1-\alpha}(c_2\alpha^2) - 3\alpha^2$ or when $c_2 = \alpha(2-\alpha)$ ($C\left(\alpha,\frac{c_2}{2}\right)$ cannot be differentiated in this point). 
Thus, for any fixed $c_2$, $c_1^*(c_2)$ can be computed by solving the equations
$$
\frac{d}{d\alpha}\left(\frac{\alpha}{1-\alpha} \frac{c_2-\alpha^2}{2}\right) = 0
$$
and
$$
\frac{d}{d\alpha} \left(\frac{2\alpha}{1-\alpha}(c_2\alpha^2) - 3\alpha^2\right) = 0,
$$
% substituting the appropriate solutions to $C(\alpha,\frac{c_2}{2})$ and selecting the maximum, comparing it with $(1-\sqrt{1-c_2})^2$ and taking the maximum. 
substituting the appropriate solutions to $C(\alpha,\frac{c_2}{2})$, and selecting the maximum out of these solutions and $(1-\sqrt{1-c_2})^2$ for $c_2\in (0,1)$. Since both equations defining the potential maximum place of $C(\alpha,\frac{c_2}{2})$ are cubic equations, this computation can be done in polynomial time with the number of digits of $c_2$ with the same precision as $c_2$ is given. We present a plot of $c_1^*$ as a function of $c_2$ in Figure \ref{fig:c1star-c2-plot}, where we also indicate the symmetric (and empirically widest) bounds.

%%% NB: (Arya) The old stuff from Section 3 that has now been replaced (e.g. my proofs for technical lemma, Angelina's initial proofs for some of the later always graphicality claims) has now been moved to the "old_proofs_and_outlines" files, for reference in case we need to pull something from there. 

%%%%%%%%%%%%%%%%%%%%%%%%%%%%%%%%%%%%%%%%%%%%%%%%%%%%%%%%%%
%%%%%%%%%%%%%%%%% § Dichotomy: NP-COMPLETE %%%%%%%%%%%%%%%
%\pagebreak
%\section{Dichotomy theorem: NP-completeness}
\section{The NP-completeness part of the dichotomy theorem}
\label{sec:np-complete}

We now proceed to the second half of our Dichotomy Theorem: namely, the NP-completeness result. 
We prove the NP-completeness by reducing the general 3-uniform hypergraphicality problem to the parameterized  3-uniform hypergraphicality problem with linear bounds $((c_1^*(c_2)-\varepsilon)){n-1\choose 2}, c_2{n-1\choose 2}$. The reduction is based on an embedding construction resembling the so-called \emph{Tyshkevich product} \cite{tyshkevich}. The first example of such embedding in the scientific literature can be found in \cite{jms}. In a nutshell, a Tyshkevich product takes two degree sequences of simple graphs, $D_1$ and $D_2$ and creates a new one $\tilde{D} := D_1\circ D_2$. The property of $\tilde{D}$ is that the number of simple graph realizations of $\tilde{D}$ is the product of the number of realizations of $D_1$ and $D_2$ \cite {barrus}. Particularly, $\tilde{D}$ is graphic if and only if both $D_1$ and $D_2$ are graphic. 

\vspace{0.5em}
It seems that the Tyshkevich product cannot be extended to 3-uniform hypergraphs in general. A heuristic explanation for this is the following. If $\tilde{D} = D_1 \circ D_2$ is a Tyshkevich product, then the realizations of $\tilde{D}$ are very rigid, and in fact, can be obtained as some (here not detailed) product of the realizations of $D_1$ and $D_2$. The proof that these are the all possible realizations of $\tilde{D}$ comes from the fact that for any degree sequence $D$, the simple graph realizations of $D$ can be transformed into each other by switch operations. The rigid structure of the realizations of $\tilde{D}$ provides that no switch operation can break this structure. Therefore, all realizations of $\tilde{D}$ have this rigid structure. On the other hand, switches are not enough to transform $3$-uniform hypergraph realizations of a degree sequence into each other \cite{hubaietal}.

\vspace{0.5em}
However, we are able to create a similar construction that builds a degree sequence $D_B$ from a general degree sequence $D_0$ and a (possibly slightly modified) degree sequence of a critical hypergraph such that $D_0$ is graphic if and only if $D_B$ is graphic. The proof that the realizations of $D_B$ are rigid is based on the rigidity of the critical hypergraphs and is significantly more involved than the proof of the rigidity of the realizations of a Tyshkevich product.

In order to formally state Part (III) of the theorem, we define the following degree sequence class, and then restate Part (II) and (III) below. 

\begin{definition}
 Given $0<c_1, c_2 <1$ and $\varepsilon > 0$,
% , and $\lambda > 0$,
 let $\generaldcce$ denote the class of $3$-uniform hypergraph degree sequences  such that for each degree sequence $D \in \generaldcce$ of length $n$, the following holds: 
\vspace{-0.5em}
 \begin{enumerate}
     \item $\sum_{d \in D} d \equiv 0 \pmod{3}$
     \item $\forall d_i \in D$, $c_1 {n-1 \choose 2} - %\lambda 
     n^{1+\varepsilon} \leq d_i \leq c_2 {n-1 \choose 2}$
\end{enumerate}
\end{definition}

\begin{theorem}\label{thm:dich-3} {\normalfont [Part (III) of \hyperref[thm:dichotomy]{Dichotomy Theorem}]}
For any $0 < c_2 < 1$ and $\varepsilon > 0$, %there exists a $\lambda > 0$ such that 
%\textit{unique}, \textit{polynomial-time computable} value $c_1^*$ ($0 < c_1^* < c_2$) such that $\forall \varepsilon > 0$,
the decision problem $\grprob$ over the degree sequence class  $\dcce$ is NP-complete, where $c_1^* := c_1^*(c_2)$ is as defined in Definition~\ref{def:c1star}.
\end{theorem}

\begin{corollary}\label{thm:dich-2} {\normalfont [Part (II) of \hyperref[thm:dichotomy]{Dichotomy Theorem}]}
For any $0 < c_2 < 1$,
%there exists a \textit{unique}, \textit{polynomial-time computable} value $c_1^*$ ($0 < c_1^* < c_2$) such that 
$\forall c_1 < c_1^*(c_2)$, $\grparam$ is NP-complete.
\end{corollary}

\begin{proof}[Proof of Theorem \ref{thm:dich-3}] 

\vspace{0.5em}
Fix $0 < c_2 < 1$ and let $c_1^*$ denote the corresponding critical value (Definition \ref{def:c1star}). We prove the theorem for %any $\lambda > 0$ and
any $0 < \epsilon < 1$ (and observe that this will imply the claim for larger $\epsilon > 1$). 

\vspace{0.5em}
%Consider any length $n \geq n_0$ for $n_0 \in \N$ that we will soon determine explicitly. Let $\dmax = \floor{c_2{n-1 \choose 2}}$ and let $k^* \in \arg \max_k f(n, k, \dmax)$. That is, $f(n, k^*, \dmax) = f^*(n, \dmax)$. Then, by Lemma \ref{lem:o-n-to-1-plus-epsilon}, there exists $n_1 \in \N$ such that $\forall n \geq n_1$, $f^*(n, \dmax) \geq c_1^* {n-1 \choose 2} - \frac{\lambda}{3}\cdot n^{1+\epsilon}$. Further, since $k^* = \Theta(n)$ by \textcolor{red}{[\textbf{TBD}: proved separately somewhere?]}, let $n_2 \in N$ be sufficiently large such that $\forall n \geq n_2$, the following holds: we can pick $\gamma > 0$ such that $\gamma^2 + \gamma < \frac{\lambda}{3}$ and $m := \floor{\gamma n^{\epsilon}}$ is less than $k$. Now we set $n_0 = \max{\{n_1, n_2, \frac{6}{\lambda}\}}$.

\vspace{0.5em}
Fixing $0<\epsilon<1$, we consider the following two decision problems. 

\vspace{0.5em}
\ul{Problem 1.} \\
{\it Input:} an arbitrary degree sequence $D \in \dcce$.\\
{\it Output:} ``Yes" if $D$ has a 3-uniform hypergraph realization, ``No" otherwise.
%$D$ of length $n$ such that $D$ is in $\dcce$, decide whether $D$ is graphic. 

\vspace{0.5em}
\ul{Problem 2.}\\
{\it Input:} an arbitrary degree sequence $D_0$.\\
{\it Output:} ``Yes" if $D_0$ has a 3-uniform hypergraph realization, ``No" otherwise.

%$D_0$ of length $m$ such that $\sum_{d \in D_0} d \equiv 0 \pmod{3}$, decide whether $D_0$ is graphic.

\vspace{0.5em}
\textbf{Reduction.} Our goal is to prove that Problem 1 is NP-complete. Clearly Problem 1 is in NP, since one can easily compute the degrees of a $3$-uniform hypergraph on $n$ vertices in polynomial time. By Deza et.\ al's result \cite{Dezaetal2019}, Problem 2 is NP-complete. % when viewed over input size $m$, but we can also view it as over input size $n$ since $m = \Theta(n^\epsilon)$ and thus $\poly(m) = \poly(n)$ \textcolor{blue}{(check reasoning)}. Thus to prove that Problem 1 is NP-complete.
We are going to prove that Problem 1 is NP-complete by showing that Problem 2 is polynomial reducible to Problem 1. In particular, we show that for any $D_0$ of length $m$, if $m$ is sufficiently large, then there exists a corresponding sequence $D_B(D_0) \in \dcce$ of length $n =  \left\lceil{}2^{\frac{1}{\varepsilon}}m^{\frac{1}{\varepsilon}}\right\rceil$ computable from $D_0$ in polynomial time such that $D_B(D_0)$ is graphic if and only if $D_0$ is graphic.

%, we show that Problem 2 is reducible to Problem 1.  In particular, we show that for any $D_0$ of length $m$, there exists a corresponding sequence $D_B(D_0) \in \dcce$ of length $n$ computable from $D_0$ in polynomial time such that $D_B(D_0)$ is graphic if and only if $D_0$ is graphic. Then, to decide Problem 2 for an arbitrary input $D_0$, we can simply compute $D_B(D_0)$ and call Problem 1 to decide graphicality. 

\vspace{0.5em}
\textbf{Critical hypergraph.} Let $D_0$ be an arbitrary degree sequence of length $m$. We first construct a critical hypergraph $H = (V_S \sqcup V_L, E)$ on $n = \left\lceil{}2^{\frac{1}{\varepsilon}}m^{\frac{1}{\varepsilon}}\right\rceil$ vertices for our chosen % $\dmax = \floor{c_2 n^2}$ 
$\dmax = \floor{c_2{n-1 \choose 2}}$ and $k = k^*(n,\dmax)$ (according to the construction described in the proof of Lemma \ref{lem:critical-degree}).
% Here  $\sqcup$ denotes the disjoint union of two vertex sets.
In particular, $|V_L| = k^*$ and $|V_S|=n-k^*$. Recall that while the construction has three cases based on $k$ and $\dmax$, we previously observed by Lemma~\ref{lem:theta-k} that for $k=k^*$, ${k-1\choose 2} < \dmax$.
%the construction is never realized in the first case (``Case 0''). 
Thus we have two cases that determine the types of edges added: 
\begin{itemize}
    \item Case 1: ${k^* - 1 \choose 2} < \dmax \leq {k^* - 1 \choose 2} + (n-k^*)(k^*-1)$
    \item Case 2: $\dmax > {k^* - 1 \choose 2} + (n-k^*)(k^*-1)$
\end{itemize}

We make one caveat in our construction of $H$ here: in Case 1 of constructing a critical hypergraph, there can in general be one 1L2S edge added, but we never add such an edge when constructing $H$. 

\vspace{0.5em}
As a result, the set of edges added in $H$ will be either (Case 1) \{all 3L, some positive number of 2L1S\}  or (Case 2) \{all 3L, all 2L1S, some positive number of 1L2S\}. This implies that $H[V_S] = \emptyset$, $H[V_L] = K_{k^*}$, and $H[V_S, V_L]$ is either (1) \{some 2L1S edges\} or (2) \{all 2L1S, some 1L2S edges\}. 

\vspace{0.5em}
Furthermore, all vertices in $V_L$ have degree $\dmax$ except for one which may have degree $\dmax - 1$ (due to the caveat above). It also follows from Definition \ref{def:critical-degree-class}, Definition \ref{def:f}, and the proof of Lemma \ref{lem:critical-degree} that for all $v \in V_S$, $d(v)$ is either $\ceil{f_0(n, k^*, \dmax)}$, $\floor{f_0(n, k^*, \dmax)}$, or possibly $\floor{f_0(n, k^*, \dmax)} - 1$ (again due to the caveat above). In particular, the minimum degree $\dmin$ in $H$ satisfies $\dmin \geq f^*(n, \dmax) - 2$. Denote the degree sequence on $n$ vertices corresponding to $H$ by $D_A$.

\vspace{0.5em}
\textbf{Embedding construction.} Let $V_N$ be an arbitrary subset of $V_L$ such that $|V_N| = m$. We now construct a hypergraph $H'$ from $H$ by first removing edges to attain $H'[V_N] = \emptyset$ and $H'[V_N, V_S] = \emptyset$. In particular, extending our edge type notation from before, this means removing all 3N, 2N1S, and 1N2S edges that were present in $H$. Then, if the $H$ construction was in Case 1 (based on $k^*$), we also remove all 1L1N1S edges present, which are a subset of the original 2L1S edges of $H$. If in Case 2, we instead keep all 1L1N1S edges, recalling that all 2L1S edges are originally present in $H$ in Case 2. When considering $H'$, we henceforth use $V_L$ to refer to the original vertex set $V_L$ of $H$ \textit{minus} the vertices $V_N$. Thus $V(H')$ is the disjoint union $V_L \sqcup V_N \sqcup V_S$. In this notation, note that $H'[V_L, V_N]$ is a complete bipartite graph, i.e. all 2L1N and 1L2N edges are still present. Denote the resulting degree sequence of $H'$ by $D_A'$. 

\vspace{0.5em}
\textbf{Constructing $D_B(D_0)$}. We now define $D_B = D_B(D_0)$ as the degree sequence obtained by adding the input sequence $D_0$ of length $m$ to the $V_N$ section of $D_A'$, recalling that $|V_N| = m$. We now claim that indeed $D_B \in \dcce$ for the 
% $\lambda$ and 
 $\epsilon$ fixed at the start. 

\begin{enumerate}
    \item First, $\sum_{d \in D_B} d \equiv 0 \pmod{3}$. This holds since $\sum_{d \in D_0} d \equiv 0 \pmod{3}$ and $\sum_{d \in D_A'} d \equiv 0 \pmod{3}$ because $D_A'$ is the degree sequence of $H'$, and by construction of $D_B$ we have $\sum_{d \in D_B} d =  (\sum_{d \in D_A'} d) + (\sum_{d \in D_0} d)$. 
    
    \vspace{0.5em}
    \item  Next, observe that the degree of each vertex only possibly decrease from $D_A$ to $D_B$. $D_A'$ is obtained from $D_A$ by removing edges. While $D_0$ is then added to the $V_N$ segment of $D_A'$ to obtain $D_B$, the magnitude of each degree's increase is at most the decrease from $D_A$ to $D_A'$ due to removing the clique $H[V_N]$. Hence each degree can only possibly decrease overall from $D_A$ to $D_B$. Since the maximum degree in $D_A$ is $\dmax = \floor{c_2{n-1 \choose 2}}$, each degree in $D_B$ is at most $c_2{n-1 \choose 2}$.

    \vspace{0.5em}
    \item Finally we prove that that $\forall d \in D_A'$, $d \geq c_1 {n-1 \choose 2} - n^{1+\epsilon}$ if $m$ (thus $n$) is sufficiently large. Since degrees only increase between $D_A'$ and $D_B$, this suffices to show that each degree in $D_B$ stays within the lower bound of $\dcce$. 
    
    \vspace{0.5em} 
    Consider the removal of edges incident to an arbitrary $v_s \in V_S$ between $D_A$ (hypergraph $H$) and $D_A'$ (hypergraph $H'$). Since we remove 2N1S, 1N2S, and possibly 1L1N1S edges, the number of edges removed incident to any $v_s$ is at most  $${|V_N| \choose 2} + |V_S| \cdot |V_N| + |V_L| \cdot |V_N| = {m \choose 2} + m n \le \left(\frac{n^{\varepsilon}}{2}\right)^2 + \frac{1}{2}n^{1+\varepsilon}
    %\leq \gamma^2 n^{2\epsilon} + \gamma n^{1+\epsilon} \leq \frac{\lambda}{3} n^{1+\epsilon},
    $$  
    %by definition of $m$ and $\gamma$. 
    which is smaller than $\frac{3}{4}n^{1+\varepsilon}$ if $m$  (thus $n$) is sufficiently large (here we use the condition that $\varepsilon < 1$).
    So $d_{H'}(v_s) \geq d_H(v_s) - \frac{3}{4} n^{1+\epsilon}$. Then, $d_{H}(v_s) \geq f^*(n, \dmax) - 2$ and $f^*(n, \floor{c_2n^2}) \geq c_1^* {n-1 \choose 2} - \frac{1}{6} n^{1+\epsilon}$ by Lemma \ref{lem:o-n-to-1-plus-epsilon}, given that $m$ (thus $n$) is sufficiently large. Thus $d_{H'}(v_s) \geq c_1^* {n-1 \choose 2} - 2 - \frac{5}{6} n^{1+\epsilon} \geq c_1^* {n-1 \choose 2} - n^{1+\epsilon}$ if $m$ (thus $n$) is sufficiently large. 

    \vspace{0.5em}
    Further, it is easy to check, by our construction of $H'$ in both Case 1 and Case 2, that for any vertex $v \notin V_S$ (i.e. $v \in V_N$ or $v \in V_L$), the vertex pairs to which $v$ is adjacent in $H'$ will be a superset of the pairs to which each $v_s \in V_S$ is adjacent. Thus, we can conclude that in fact for all $v \in V(H')$, $d_{H'}(v) \geq c_1^* {n-1 \choose 2} - n^{1+\epsilon}$, and then we are done since $D_A'$ is the degree sequence of $H'$. 
\end{enumerate}

Thus we have shown that $D_B \in \dcce$. Also, the reduction can be done in polynomial time since $\varepsilon$ is a fixed positive constant.

\vspace{0.5em}
{\bf The reduction.} Next we are going to prove that $D_B$ is graphic if and only if $D_0$ is graphic.
%It now just remains to show that $D_B$  is graphic if and only if $D_0$ is graphic. This is proved subsequently in Lemma \ref{lem:DBD0}. \end{proof}
Before we prove it, we introduce some notation. If $V = V_S\sqcup V_N \sqcup V_L$ (that is, disjoint union of small, intermediate, and large degree vertices), then a degree sequence $D$ on it can be split into disjoint union of sequences $D[S]$, $D[N]$, and $D[L]$. Regarding the three types of vertices a hyperedge can be incident with, there are ${3+3-1\choose 3} = 10$ different types of hyperedges that we will denote by 3L, 2L1N, etc. similarly as in Definition~\ref{def:hyperegde-types}.  We will denote the total degree sum due to an edge type ``$ABC$'' on vertex part ``$V_X$'' in hypergraph ``$F$'' as $[ABC]^F_X$. For instance, $[\text{2L1S}]^{H'}_L$ denotes the total degree sum on $V_L$ contributed by 2L1S edges in $H'$ and is equal to twice the number of 2L1S edges. %, while $[\text{2L1S}]^{H'}_N = 0$.

 %We will refer to the $S$, $N$, and $L$ parts of the degree sequence $D_B$, denoted $D_B[S]$, $D_B[N]$, and $D_B[L]$, by which we mean the sequence of degrees corresponding to $V_S$, $V_N$, and $V_L$ vertex parts of $H'$. In particular $D_B$ is the concatenation of $D_B[S]$, $D_B[N]$, and $D_B[L]$. We will also refer to the $V_S$, $V_N$, $V_L$ components of $G$, by which we mean the set of vertices corresponding to the degree sequence parts $D_B[S]$, $D_B[N]$, and $D_B[L]$ respectively. Hence we use $V_S, V_N, V_L$ refer to the corresponding vertex parts of both $H'$ and $G$; the hypergraph being referred to will be clear by context.   

%\begin{lemma}\label{lem:DBD0}
%The degree sequence $D_B = D_B(D_0)$ on $n$ vertices is graphic if and only if the sequence $D_0$ on $m$ vertices is graphic. 
%\end{lemma}
% Notation (Arya): if someone has a better idea for notation for the different degree sequences / hypergraphs used in the course of the logic for NP-completeness please feel free to comment. This was the most consistent thing I could come up with from thinking for a bit but it's not great 

%\begin{lemmaproof}
%\begin{proof}
One direction ($\Leftarrow$) is trivial: if $D_0$ is graphic, then by the given construction, 
%outlined in the proof of Theorem \ref{thm:dich-3}, 
$D_B$ must be graphic. Specifically, letting $G_0$ be a realization of $D_0$ on the vertices $V_N$, the hypergraph $H' \sqcup G_0$ is a realization of $D_B$. 
\vspace{0.5em}

For the other direction ($\Rightarrow$), suppose $D_B$ is graphic, and let $G$ be a realization of $D_B$. %\textcolor{orange}{(move this up to main proof? $\rightarrow$)} We will refer to the $S$, $N$, and $L$ parts of the degree sequence $D_B$, denoted $D_B[S]$, $D_B[N]$, and $D_B[L]$, by which we mean the sequence of degrees corresponding to $V_S$, $V_N$, and $V_L$ vertex parts of $H'$. In particular $D_B$ is the concatenation of $D_B[S]$, $D_B[N]$, and $D_B[L]$. We will also refer to the $V_S$, $V_N$, $V_L$ components of $G$, by which we mean the set of vertices corresponding to the degree sequence parts $D_B[S]$, $D_B[N]$, and $D_B[L]$ respectively. Hence we use $V_S, V_N, V_L$ refer to the corresponding vertex parts of both $H'$ and $G$; the hypergraph being referred to will be clear by context.   
We will prove the claim by showing that the degree sequence of $G[V_N]$ is $D_0$.
Considering the hypergraph $G$, let $$\Sigma_S = \sum_{d_i \in D_B[S]} d_i = \sum_{v \in V_S} d(v) \text{ and } \Sigma_L = \sum_{d_i \in D_B[L]} d_i = \sum_{v \in V_L} d(v).$$ 

Let $G'$ denote the hypergraph obtained by removing all 3N edges from $G$, i.e. $G'[V_N] = \emptyset$, and let $D_B'$ denote the degree sequence of $G'$. $\Sigma_S$ and $\Sigma_L$ stay the same in $G'$ since the degrees of vertices  $V_S$ and $V_L$ do not change when removing 3N edges. Furthermore, recall from our previous definitions that the hypergraph $H'$ is a realization of the degree sequence $D_A'$, where adding $D_0$ to the $V_N$ section of $D_A'$ yields the sequence $D_B$. The $V_L$ and $V_S$ sections are not modified by the addition of $D_0$, and so $\Sigma_S = \sum_{d_i \in D_A'[S]} d_i$ and $\Sigma_L = \sum_{d_i \in D_A'[L]} d_i$. 
\vspace{0.5em}

We now analyze the ratio $\frac{\Sigma_L}{\Sigma_S}$ according to the types of edges present in two cases, corresponding to Case 1 and Case 2 in the construction of $H'$. Since $H'[N] = \emptyset$ and $G'[N] = \emptyset$, there are no 3N edges in either, so there are 9 possible edge types: 3L, 2L1N, 1L2N, 2L1S, 1L1N1S, 1S2L, 2N1S, 1N2S, 3S.% We will denote the total degree sum due to an edge type ``$ABC$'' on vertex part ``$V_X$'' in hypergraph ``$F$'' as $[ABC]^F_X$. For instance, $[\text{2L1S}]^{H'}_L$ denotes the total degree sum on $V_L$ contributed by 2L1S edges in $H'$ and is equal to twice the number of 2L1S edges, while $[\text{2L1S}]^{H'}_N = 0$.
\vspace{0.5em}

In each case, we will show that $D_B'[N] = D_A'[N]$, and furthermore, this degree sequence segment is regular. This will suffice for the proof because  we know that adding $D_0$ to the $V_N$ segment of $D_A'$ (i.e. $D_A'[N]$) yields $D_B$ overall. Given the regularity of $D_B'[N] = D_A'[N]$, this implies that $D_B[N] - D_B'[N] = D_0$, while we also know that $D_B[N] - D_B'[N]$ must be precisely the degree sequence of the induced subgraph $G[V_N]$ removed from $G$ to form $G'$. That is, the subgraph removed from $G$ is a realization of $D_0$, so $D_0$ is graphic if $D_B$ is graphic.  

\vspace{0.5em}
\textbf{Case 1.} The edge set of $H'$ is given by all possible 3L, 2L1N, and 1L2N edges; \textit{some} 2L1S edges according to the critical hypergraph construction; and no 1L1N1S, 1L2S, 2S1N, 1S2N, 3S, or 3N edges.\footnote{Recall that while the critical hypergraph in Case 1 might have one 1L2S edge for rounding, the construction of $H'$ importantly drops this edge.} In particular

$$ \frac{\Sigma_L}{\Sigma_S} = \frac{\sum_{d_i \in D_A'[L]} d_i}{\sum_{d_i \in D_A'[S]} d_i} = \frac{[\text{3L}]^{H'}_L + [\text{2L1N}]^{H'}_L  + [\text{1L2N}]^{H'}_L + [\text{2L1S}]^{H'}_L}{[\text{2L1S}]^{H'}_S} > 2  $$

$$ \text{and } \frac{\Sigma_L - [\text{3L}]^{H'}_L - [\text{2L1N}]^{H'}_L  - [\text{1L2N}]^{H'}_L}{\Sigma_S} = \frac{[\text{2L1S}]^{H'}_L}{[\text{2L1S}]^{H'}_S} = 2$$

This holds because each edge of these types adds either at least $2$ degrees to $\Sigma_L$ while adding at most $1$ degree to $\Sigma_S$, or adds at least $1$ degree to $\Sigma_L$ while adding $0$ degrees to $\Sigma_S$. The inequality is strict because there is a positive number of 3L edges (and 2L1N edges, and 1L2N edges, indeed). Further, $[\text{2L1S}]^{H'}_L = 2 \cdot [\text{2L1S}]^{H'}_S$ and this value is positive (non-zero) because we are in Case 1 of $H'$. 
\vspace{0.5em}

Recall that the sizes of vertex sets $|V_S|$, $|V_N|$, $|V_L|$ are the same in $H'$ and $G'$, and so there are the same number of edges of each type, and $H'$ has all possible 3L, 2L1N, and 1L2N edges. It follows that $[\text{3L}]^{G'}_L \leq [\text{3L}]^{H'}_L$, $ [\text{2L1N}]^{G'}_L \leq [\text{2L1N}]^{H'}_L$, and $[\text{1L2N}]^{G'}_L \leq [\text{1L2N}]^{H'}_L$. But then 
$$ \frac{\Sigma_L - [\text{3L}]^{G'}_L - [\text{2L1N}]^{G'}_L  - [\text{1L2N}]^{G'}_L}{\Sigma_S} \geq \frac{\Sigma_L - [\text{3L}]^{H'}_L - [\text{2L1N}]^{H'}_L  - [\text{1L2N}]^{H'}_L}{\Sigma_S} = 2 \ \ (\mathbf{\star})$$
with equality if and only if $[\text{3L}]^{G'}_L = [\text{3L}]^{H'}_L$, $ [\text{2L1N}]^{G'}_L = [\text{2L1N}]^{H'}_L$, and $[\text{1L2N}]^{G'}_L = [\text{1L2N}]^{H'}_L$. 
\vspace{0.5em}

We will show that equality must hold. Aside from 3L, 2L1N, and 1L2N edges, the remaining edge types in $G'$ are $\mathcal{E} = \{\text{2L1S, 1L1S1N, 1L2S, 2N1S, 1N2S, 3S}\}$. Considering $\Sigma_L$ and $\Sigma_S$ in the context of $G'$, the numerator of $\frac{\Sigma_L - [\text{3L}]^{G'}_L - [\text{2L1N}]^{G'}_L  - [\text{1L2N}]^{G'}_L}{\Sigma_S}$ is the total degree contributed by edges of types from $\mathcal{E}$ to $V_L$, and the denominator is the total degree contributed to $V_S$ (from edges of these types). Notice that an edge of type in $\mathcal{E}$ contributes \textit{at most} a ratio of $2$ in terms of (degrees contributed to $V_L$)/(degrees contributed to $V_S$). It follows that $(\Sigma_L - [\text{3L}]^{G'}_L - [\text{2L1N}]^{G'}_L  - [\text{1L2N}]^{G'}_L)/(\Sigma_S) \leq 2$. Combining this with $(\star)$ shows that equality must hold. 
\vspace{0.5em}

Hence $[\text{3L}]^{G'}_L = [\text{3L}]^{H'}_L$, $ [\text{2L1N}]^{G'}_L = [\text{2L1N}]^{H'}_L$, $[\text{1L2N}]^{G'}_L = [\text{1L2N}]^{H'}_L$. Furthermore, 2L1S edges are the only type in $\mathcal{E}$ with a ratio of $2$ for large degrees to small degrees. As such, for equality to hold it must be true that the only edges in $G'$ aside from 3L, 2L1N, 1L2N are of type 2L1S. The crucial point is that the edges incident with $V_N$ in $G'$ are exactly all possible 2L1N and 1L2N edges, as is the case in $H'$. Hence $D_B'[N] = D_A'[N]$, and this degree sequence segment is clearly regular since the complete subgraphs of the two types are present. 

\vspace{0.5em} 
\textbf{Case 2.} The argument is largely analogous to Case 1. The edge set of $H'$ is given by all possible 3L, 2L1N, 1L2N, 2L1S, 1L1N1S edges; \textit{some} 1L2S edges according to the critical hypergraph construction; and no 2S1N, 1S2N, 3S, or 3N edges. In particular 
$$ \frac{\Sigma_L}{\Sigma_S} = \frac{[\text{3L}]^{H'}_L + [\text{2L1N}]^{H'}_L  + [\text{1L2N}]^{H'}_L + [\text{2L1S}]^{H'}_L + [\text{1L1N1S}]^{H'}_L + [\text{1L2S}]^{H'}_L}{[\text{2L1S}]^{H'}_S+ [\text{1L1N1S}]^{H'}_S + [\text{1L2S}]^{H'}_S} > \frac{1}{2}  $$
$$ \text{and } \frac{\Sigma_L - [\text{3L}]^{H'}_L - [\text{2L1N}]^{H'}_L  - [\text{1L2N}]^{H'}_L - [\text{1L1N1S}]^{H'}_L}{\Sigma_S  - [\text{2L1S}]^{H'}_S - [\text{1L1N1S}]^{H'}_S } = \frac{[\text{1L2S}]^{H'}_L}{[\text{1L2S}]^{H'}_S} = \frac{1}{2}.$$

This holds because each edge of these types adds either at least $2$ degrees to $\Sigma_S$ while adding at most $1$ degree to $\Sigma_L$, or adds at least $1$ degree to $\Sigma_S$ while adding $0$ degrees to $\Sigma_L$. Furthermore, there is a positive number of 2S1L edges since we are in Case 2 of $H'$. 
\vspace{0.5em} 

Since $H'$ has all possible 3L, 2L1N, 1L2N, 2L1S, and 1L1N1S edges, we have that $[\text{3L}]^{G'}_L \leq [\text{3L}]^{H'}_L, [\text{2L1N}]^{G'}_L \leq [\text{2L1N}]^{H'}_L, [\text{1L2N}]^{G'}_L \leq [\text{1L2N}]^{H'}_L, [\text{2L1S}]^{G'}_L \leq [\text{2L1S}]^{H'}_L, [\text{1L1N1S}]^{G'}_L \leq [\text{1L1N1S}]^{H'}_L$ and so $$\frac{\Sigma_L - [\text{3L}]^{G'}_L - [\text{2L1N}]^{G'}_L  - [\text{1L2N}]^{G'}_L - [\text{1L1N1S}]^{G'}_L}{\Sigma_S  - [\text{2L1S}]^{G'}_S - [\text{1L1N1S}]^{G'}_S} \geq \frac{1}{2} \ \ ( \star \star) $$ with equality if and only if equality holds in all of the prior inequalities. Indeed equality holds in $(\star \star)$ due to an exactly analogous argument to the one in Case 1. The remaining edge types in this case are $\mathcal{E} = \{\text{1L2S, 2N1S, 1N2S, 3S}\}$. These types have \textit{at most} a ratio of $\frac{1}{2}$ for large degrees to small degrees contributed, and so $(\Sigma_L - [\text{3L}]^{G'}_L - [\text{2L1N}]^{G'}_L  - [\text{1L2N}]^{G'}_L - [\text{1L1N1S}]^{G'}_L)/(\Sigma_S  - [\text{2L1S}]^{G'}_S - [\text{1L1N1S}]^{G'}_S) \leq \frac{1}{2}$. 
\vspace{0.5em}

Thus $[\text{3L}]^{G'}_L = [\text{3L}]^{H'}_L, [\text{2L1N}]^{G'}_L = [\text{2L1N}]^{H'}_L, [\text{1L2N}]^{G'}_L = [\text{1L2N}]^{H'}_L, [\text{2L1S}]^{G'}_L = [\text{2L1S}]^{H'}_L, [\text{1L1N1S}]^{G'}_L = [\text{1L1N1S}]^{H'}_L$. Furthermore, for equality to hold, the only edges in $G'$ aside from these types must be of type 1L2S, since 1L2S edges are the only type in $\mathcal{E}$ with a ratio of $\frac{1}{2}$ for large degrees to small degrees. In particular, the edges adjacent to $V_N$ in $G'$ are exactly the set of all possible 2L1N, 1L2N, and 1L1N1S edges, as in $H'$. Hence $D_B'[N] = D_A'[N]$, and this degree sequence segment is again regular since the complete subgraphs of these three types are present. % \end{lemmaproof}
\end{proof}

%\pagebreak
\section{Asymptotic always graphic interval bounds in $t$-uniformity}   \label{sec:asymp}

We now briefly consider the question of characterizing always graphic intervals when we move from $3$-uniformity to general case of $t$-uniformity for arbitrary $t$. Informally, our main result shows that the width of always graphic intervals diminishes to $0$ asymptotically in $t$. Concretely, we show that for $t$-uniform hypergraphs, the width of any always graphic interval\footnote{In particular, centered at any point, rather than just symmetric intervals around $\frac{1}{2}$} is bounded by $O(t^{-\frac{1}{3}})$; this is formally stated in Theorem \ref{thm:asymp-t}. We leave attempting to prove a complete dichotomy theorem for $t$-uniform hypergraphs for future work. 

\begin{theorem}\label{thm:asymp-t}
There exists a function $c(t) = C \cdot t^{-\frac{1}{3}}$ for a constant $C$ such that for any $t > 1$, for any center $p$ and sufficiently large $n$, there exists a \textit{non-graphic} $t$-uniform hypergraphic degree sequence of length $n$ with all degrees between $[(p-c(t)) \cdot {n-1 \choose t-1}, (p+c(t)) \cdot {n-1 \choose t-1}]$. 
\end{theorem}

% how this leads to corollary + add discussion on how theorem c(t) is not particularly optimized/minimized + why we don't pursue a dichotomy theorem for general t. 
Corollary \ref{cor:asympt-t} below concretely states the main high-level result of this section. Observe that this follows immediately from Theorem \ref{thm:asymp-t} since $c(t) = C \cdot t^{-\frac{1}{3}}$ tends to $0$ as $t \rightarrow \infty$, and the theorem shows that there can be no always graphic interval in $t$-uniformity with width greater than $2 \cdot c(t)$. We remark that the analysis to derive the specific $c(t)$ function used in our results is not particularly optimized, as even this $c(t)$ function obtained through basic analysis is sufficient to conclude Corollary \ref{cor:asympt-t}. 

\begin{corollary}\label{cor:asympt-t} 
The width of the largest always graphic interval for $t$-uniform hypergrapic degree sequences goes to $0$ as $t \rightarrow \infty$. 
% Update (Arya): relative to old corollary, just removing the "symmetric" part
\end{corollary}

The proof of Theorem \ref{thm:asymp-t} uses the technical result presented in Lemma \ref{lem:asymp-t-prob} below. The proof of the lemma uses a straightforward concentration argument applied to a hypergeometric distribution and is deferred to \ref{sec:appendix-proofs}. 

\begin{lemma}\label{lem:asymp-t-prob}
Fix arbitrary $t\in \N$ such that $t > 1$. Let $\epsilon > 0$ be arbitrary, and consider any $\delta > 0$ such that $\delta \geq (\epsilon \cdot t)^{-\frac{1}{2}}$. 
% let $\delta = \left(\frac{t}{2}\right)^{-\frac{1}{3}}$ and $\epsilon = (4t)^{-\frac{1}{3}}$. 
Then for any even $n \geq t$, the following holds: if $V := A \sqcup B$ is a ground set of size $n$ with $|A| = |B| = \frac{n}{2}$, and $W$ is a $t$-subset drawn uniformly at random from $V$, then 
$\pr{||W \cap A| - |W \cap B|| > \delta t} < \epsilon$. 
% Update (Arya): the \epsilon and \delta values here are different from before (see gray outline below) -- constant factors change a bit with the simpler proof directly from hypergeometric + chebyshev. Also, removed the "there exists $n_0 \geq t$ such that for any even $n \geq n_0$" -- now don't need "large enough n" claim?
\end{lemma}

\begin{proof}[Proof of Theorem \ref{thm:asymp-t}]

\vspace{0.5em}
Let $C_0 = (2^{\frac{1}{3}} + 4^{-\frac{1}{3}})$, and let $C$ be any constant such that $C > C_0$. Consider arbitrary $t > 1$. Let $\delta = \left(\frac{t}{2}\right)^{-\frac{1}{3}}$ and $\epsilon = (4t)^{-\frac{1}{3}}$, and define $c(t) = C \cdot t^{-\frac{1}{3}}$. 

\vspace{0.5em}
First, since $C > C_0$, observe that we can consider $n$ sufficiently large such that (a) $n \geq t$ and (b) $\floor{(p+c(t)) \cdot {n-1 \choose t-1}} - \ceil{(p-c(t)) \cdot {n-1 \choose t-1}} > 2  C_0 \cdot t^{-\frac{1}{3}} \cdot {n-1 \choose t-1} = 2(\epsilon + \delta)\cdot {n-1 \choose t-1}$. Thus define $\dmin = \ceil{(p-c(t)) \cdot {n-1 \choose t-1}}$ and $\dmax = \floor{(p+c(t)) \cdot {n-1 \choose t-1}}$, such that $\dmax - \dmin > 2(\epsilon + \delta)\cdot {n-1 \choose t-1}$. Furthermore, we assume WLOG that $n$ is even. We will prove that the $n$-length sequence $D$ given by $\frac{n}{2}$ $\dmin$ degrees and $\frac{n}{2}$ $\dmax$ degrees is not graphic, which implies the claim since $\dmin, \dmax \in [(p-c(t)) \cdot {n-1 \choose t-1}, (p+c(t)) \cdot {n-1 \choose t-1}]$.

\vspace{0.5em}
Assume for contradiction that there exists a $t$-uniform hypergraph $G$ on $n$ vertices $V = V_S \sqcup V_L$, $|V_S| = |V_L| = \frac{n}{2}$, such that $d(v) = \dmin$ for all $v \in V_S$ and $d(v) = \dmax$ for all $v \in V_L$. We will first translate the probability statement from Lemma \ref{lem:asymp-t-prob} into a counting statement to bound the number of possible edges that, informally, could contribute a large difference between the degrees of $V_S$ and $V_L$. 

\vspace{0.5em}
Let $E$ denote the set of edges of the \textit{complete} $t$-uniform hypergraph over the same vertex set $V = V_S \sqcup V_L$, and for any $t_1 + t_2 = t$, let $E_{t_1, t_2} \subseteq E$ be the edges that contain $t_1$ $V_L$ vertices and $t_2$ $V_S$ vertices, that is, $E_{t_1, t_2} := \{e \in E: |e \cap V_L| = t_1, |e \cap V_S| = t_2 \}$. Then let $\EH = \bigcup_{|t_1 - t_2| > \delta t} E_{t_1, t_2}$ and $\EL = \bigcup_{|t_1 - t_2| \leq \delta t} E_{t_1, t_2} = E\setminus \EH$. Suppose $W$ is an edge drawn uniformly at random from $E$, or equivalently, a $t$-subset drawn uniformly at random from $V$. Observe that the event $\{||W \cap V_{L}| - |W \cap V_{S}|| > \delta \cdot t\}$ is equal to the event $\{W \in \EH\}$. In particular, we have that $$\pr{||W \cap V_{L}| - |W \cap V_{S}|| > \delta \cdot t} = \pr{W \in \EH} = \frac{|\EH|}{|E|}.$$ Since $\delta \geq (\epsilon \cdot t)^{-\frac{1}{2}}$, by Lemma \ref{lem:asymp-t-prob} it follows that $\frac{|\EH|}{|E|} < \epsilon \Rightarrow |\EH| < \epsilon \cdot {n \choose t}$.

\vspace{0.5em}
We can now bound the maximum difference achievable between the $\dmin$ and $\dmax$ degrees. Let $\Delta := \left(\sum_{v \in V_L} d(v)\right) -  \left(\sum_{v \in V_S} d(v)\right) = \frac{n}{2} (\dmax - \dmin)$ denote the total degree difference between vertices in $V_L$ and $V_S$. To upper bound the contribution of the edges $E(G) \subseteq E = \EH \sqcup \EL$, we will separately upper bound (loosely) the possible contributions of $\EH$ edges and $\EL$ edges. Observe that an edge $e \in E_{t_1, t_2}$ contributes exactly $t_1 - t_2$ to $\Delta$. Trivially this means that any edge $e$ can contribute at most $t$ to $\Delta$, and so we can bound the contribution of $\EH$ edges to $\Delta$ by $|\EH| \cdot t < \epsilon t \cdot {n \choose t}$. Since any edge $e \in \EL$ cannot contribute more than $|t_1 - t_2| < \delta t$ to $\Delta$, we can bound the contribution by $\EL$ edges to $\Delta$ by $|\EL| \cdot \delta t \leq \delta t \cdot {n \choose t}$. Thus $\frac{n}{2} (\dmax - \dmin) = \Delta < (\epsilon + \delta) \cdot t \cdot {n \choose t}$, which by rearranging yields $\dmax - \dmin < 2 (\epsilon + \delta) {n-1 \choose t-1}$. This is a contradiction and completes the proof. \end{proof}

\section{Acknowledgements}
This project began as a part of Budapest Semesters in Mathematics in Spring 2024, and the authors thank BSM for running the program. IM was supported by NKFIH grants SNN135643 and  K132696.

%%% NB: (Arya) The old versions of asymptotic claim proofs have now been moved to the "old_proofs_and_outlines" files, for reference in case we need to pull something from there. 

%% REFERENCES %% 

\pagebreak
\appendix
\section*{\appendixname}

\section{Omitted Proofs from Section \ref{sec:always=graphic}}\label{app:3}
%\textcolor{orange}{TBD: move technical lemma proofs to here.}

\subsection{Proof of Lemma~\ref{lem:discrete-continuous-f0-g}}

\begin{proof}
First observe by Definition \ref{def:g} that the difference between $f_0\left(n,\floor{\alpha n},\floor{c_2 {n-1\choose 2}}\right)$ and $g\left(n,\floor{\alpha n},\floor{c_2 {n-1\choose 2}}\right)$ is $O(n)$. Therefore
$$
\lim_{n\rightarrow \infty} \frac{f_0\left(n,\floor{\alpha n},\floor{c_2 {n-1\choose 2}}\right) - g\left(n,\floor{\alpha n},\floor{c_2 {n-1\choose 2}}\right)}{{n-1\choose 2}} = 0.
$$
Therefore, it suffices to prove that
$$
\lim_{n\rightarrow \infty} \frac{f_0\left(n,\floor{\alpha n},\floor{c_2 {n-1\choose 2}}\right)}{{n-1\choose 2}}  = C\left(\alpha, \frac{c_2}{2}\right)
$$
and then the limit of $\frac{g\left(n,\floor{\alpha n},\floor{c_2 {n-1\choose 2}}\right)}{{n-1\choose 2}}$ follows.

From Definition \ref{def:critical-degree-class} and Definition \ref{def:f}, we can concretely write 

\begin{equation*}
f_0(n, k, \dmax) = \begin{cases} 
  \frac{1}{n-k} \left(\pmodtight{2k\dmax}{3}\right) & \text{[C0]} \\
  
  \frac{1}{n-k}\left(\left\lfloor\frac{k\left(\dmax-{k-1\choose 2}\right)}{2}\right\rfloor + 2\left(\pmodtight{k\left(\dmax-{k-1\choose 2}\right)}{2}\right)\right) & \text{[C1]} \\
  
  {k\choose 2} + \frac{1}{n-k} \left(2k\left(\dmax - {k-1\choose 2} - (n-k)(k-1)\right)\right) & \text{[C2]} 
% f_0(n, k, \dmax) = \begin{cases} 
%   0 & \dmax\leq \binom{k-1}{2}\\
%   \frac{k}{2(n-k)}\left(\dmax-\binom{k-1}{2}\right) & \binom{k-1}{2}<\dmax\leq \binom{k-1}{2}+(n-k)(k-1)\\
%   \frac{2k}{n-k}\left(\dmax-\binom{k-1}{2}-(n-k)(k-1) \right)+\binom{k}{2} & \binom{k-1}{2}+(n-k)(k-1)<\dmax
\end{cases}
\end{equation*} 

where [C0] is $\dmax\leq \binom{k-1}{2}$, [C1] is $\binom{k-1}{2}<\dmax\leq \binom{k-1}{2}+(n-k)(k-1)$, and [C2] is $\binom{k-1}{2}+(n-k)(k-1)<\dmax$.

\vspace{0.5em}

    If $\frac{c_2}{2} < \frac{\alpha^2}{2}$, then for sufficiently large $n$, $\floor{c_2 {n-1\choose 2}} \le {\floor{\alpha n}-1 \choose 2}$. Then condition [C0] holds, and thus
    $$
    f_0\left(n,\floor{\alpha n},\floor{c_2 {n-1\choose 2}}\right) = \frac{1}{n-\floor{\alpha n}} \left(\pmodtight{2\floor{\alpha n}\floor{c_2 {n-1\choose 2}}}{3}\right).
    $$
    Then clearly
    $$
    \lim_{n\rightarrow \infty} \frac{f_0\left(n,\floor{\alpha n},\floor{c_2 {n-1\choose 2}}\right)}{{n-1\choose 2}} = 0 = C\left(\alpha, \frac{c_2}{2}\right).
    $$

    If $\frac{c_2}{2} =\frac{\alpha^2}{2}$, then for some $n$, condition [C0] holds, while for other $n$'s, condition [C1] holds. For the subset of $n$'s for which condition [C0] holds, the considered limit is clearly $0$ by the same reasoning as above. 
    % Consider the subset of $n$'s for which condition [C0] holds. For that subsequence of $n$'s, the considered limit is clearly $0$. 
    Now consider the $n$'s for which condition [C1] holds. In this case, we have 
    $$
    f_0\left(n,\floor{\alpha n},\floor{c_2 {n-1\choose 2}}\right) = \frac{1}{n-\floor{\alpha n}}\left(\left\lfloor\frac{\floor{\alpha n}\left(\floor{c_2 {n-1\choose 2}}-{\floor{\alpha n}-1\choose 2}\right)}{2}\right\rfloor + \right.$$
    $$
    \left. 2\left(\pmodtight{\floor{\alpha n}\left(\floor{c_2 {n-1\choose 2}}-{\floor{\alpha n}-1\choose 2}\right)}{2}\right)\right)
    $$
    $$
    \le \frac{1}{n-\floor{\alpha n}}\left(\left\lfloor\frac{\floor{\alpha n}\left(\floor{c_2 {n-1\choose 2}}-{\floor{\alpha n}-1\choose 2}\right)}{2}\right\rfloor 
    +  2\right).
    $$
    
    Since $\floor{c_2 {n-1\choose 2}}-{\floor{\alpha n}-1\choose 2} = O(n)$ (recall that $c_2 = \alpha^2$) and $\frac{\floor{\alpha n}}{n- \floor{\alpha n}} = O(1)$,
     $$
    \lim_{n\rightarrow \infty} \frac{f_0\left(n,\floor{\alpha n},\floor{c_2 {n-1\choose 2}}\right)}{{n-1\choose 2}} = 0 = C\left(\alpha, \frac{c_2}{2}\right).
    $$

    If $\frac{\alpha^2}{2} < \frac{c_2}{2} < \alpha \left( 1-\frac{\alpha}{2}\right)$, then for sufficiently large $n$,
    $$
    {\floor{\alpha n}-1\choose 2} < \floor{c_2 {n-1\choose 2}}  < {\floor{\alpha n}-1 \choose 2} + ( n-\floor{\alpha n})(\floor{\alpha n} -1),
    $$
    and thus condition [C1] holds. Therefore,
     $$
    f_0\left(n,\floor{\alpha n},\floor{c_2 {n-1\choose 2}}\right) = \frac{1}{n-\floor{\alpha n}}\left(\left\lfloor\frac{\floor{\alpha n}\left(\floor{c_2 {n-1\choose 2}}-{\floor{\alpha n}-1\choose 2}\right)}{2}\right\rfloor + \right.$$
    $$
    \left. 2\left(\pmodtight{\floor{\alpha n}\left(\floor{c_2 {n-1\choose 2}}-{\floor{\alpha n}-1\choose 2}\right)}{2}\right)\right)
    $$
    $$
    \le \frac{1}{n-\floor{\alpha n}}\left(\left\lfloor\frac{\floor{\alpha n}\left(\floor{c_2 {n-1\choose 2}}-{\floor{\alpha n}-1\choose 2}\right)}{2}\right\rfloor 
    +  2\right).
    $$ 
    Then, by expanding the expressions it follows that 
    $$
    \lim_{n\rightarrow \infty}\frac{\frac{1}{n-\floor{\alpha n}}\left(\left\lfloor\frac{\floor{\alpha n}\left(\floor{c_2 {n-1\choose 2}}-{\floor{\alpha n}-1\choose 2}\right)}{2}\right\rfloor 
    +  2\right)}{{n-1\choose 2}} = \frac{\alpha}{1-\alpha} \left(\frac{c_2-\alpha^2}{2}\right) = C(\alpha ,\frac{c_2}{2})
    $$

    If $\frac{c_2}{2} = \alpha \left( 1-\frac{\alpha}{2}\right)$, then for some $n$, condition [C1] holds, while for other $n$'s, condition [C2] holds. % Consider the subset of $n$'s for which condition [C1] holds. For that subsequence of $n$'s, the considered limit is clearly $C(\alpha ,\frac{c_2}{2})$.
    For the subset of $n$'s for which condition [C1] holds, the considered limit is clearly $C(\alpha ,\frac{c_2}{2})$ by the same reasoning as above. Now consider the $n$'s for which condition [C2] holds. In this case, $f_0\left(n,\floor{\alpha n},\floor{c_2 {n-1\choose 2}}\right)$ is equal to the following expression:
     $$
     {\floor{\alpha n}\choose 2} +\frac{1}{n-\floor{\alpha n}} \left(2\floor{\alpha n}\left(\floor{c_2 {n-1\choose 2}} - {\floor{\alpha n}-1\choose 2} - (n-\floor{\alpha n})(\floor{\alpha n}-1)\right)\right).
    $$
    Hence we can write
    \begin{align*}
    \lim_{n\rightarrow \infty} \frac{f_0\left(n,\floor{\alpha n},\floor{c_2 {n-1\choose 2}}\right)}{{n-1\choose 2}} &= \alpha^2 +\frac{2\alpha}{1-\alpha}\left(c_2 - \alpha^2 - 2(1-\alpha)\alpha\right) \\
    &= \frac{\alpha^2 (1-\alpha) +2\alpha (\alpha(2-\alpha) - \alpha^2 - 2(1-\alpha)\alpha)}{1-\alpha} \\
    &= \frac{\alpha (\alpha -\alpha^2)}{1-\alpha} \\
    &=  \frac{\alpha}{1-\alpha} \cdot \frac{\alpha(2-\alpha) - \alpha^2}{2} \\
    &= \frac{\alpha}{1-\alpha} \cdot \frac{c_2 - \alpha^2}{2} \\
    &= C(\alpha,\frac{c_2}{2}).
    \end{align*}

    Finally, if  $\frac{c_2}{2} > \alpha \left( 1-\frac{\alpha}{2}\right)$, then for sufficiently large $n$, condition [C2] holds. Then 
    $$
     f_0\left(n,\floor{\alpha n},\floor{c_2 {n-1\choose 2}}\right) =  
    $$
    $$
    {\floor{\alpha n}\choose 2} + \frac{1}{n-\floor{\alpha n}} \left(2\floor{\alpha n}\left(\floor{c_2 {n-1\choose 2}} - {\floor{\alpha n}-1\choose 2} - (n-\floor{\alpha n})(\floor{\alpha n}-1)\right)\right).
    $$
    Then
    $$
    \lim_{n\rightarrow \infty} \frac{f_0\left(n,\floor{\alpha n},\floor{c_2 {n-1\choose 2}}\right)}{{n-1\choose 2}} 
    = \alpha^2 +\frac{2\alpha}{1-\alpha}\left(c_2 - \alpha^2 - 2(1-\alpha)\alpha\right) =
    $$
    $$
    = \frac{2\alpha}{1-\alpha}(c_2 - \alpha^2) - 3\alpha^2 = C(\alpha,\frac{c_2}{2}).
    $$

    To prove the speed of convergence, observe that in each case, both the $f_0$ and the $g$ functions can be lower and upper bounded by fractions of polynomials (of $n$, while $\alpha$ and $c_2$ are constants). That is $\alpha n - 1\le \floor{\alpha n} \le \alpha n$, etc., and the modular function parts can be lower and upper bounded by constants. The proof of convergences is not detailed in this proof, but it can be shown by the Squeeze Theorem using these fractions of polynomials. Each polynomial is an order of at most $3$, the coefficients are bounded (in fact, each coefficient is between $0$ and $1$), and the limit values are bounded ($C(\alpha,d)$ is bounded between $0$ and $1$). Then it is easy to see that fractions of polynomials converge to their limit polynomially quickly, and there is a universal polynomial for the speed of convergence.
\end{proof}

\subsection{Proof of Lemma~\ref{lem:conv-f-star-g-star}}

\begin{proof}
    We are going to prove these limits by definition. That is, we show that for any $\varepsilon > 0$, there exists an $n_0$ such that for any $n \ge n_0$,
    $$
    \left|\frac{f^*(n,\floor{c_2 {n-1\choose 2}})}{{n-1\choose 2}} - c_1^*(c_2) \right| \leq \varepsilon
    $$
    and
    $$
    \left|\frac{g^*(n,\floor{c_2 {n-1\choose 2}})}{{n-1\choose 2}} - c_1^*(c_2) \right| \leq \varepsilon
    $$
    We show the proof for the first limit.
    Fix some $c_2$ and $\varepsilon >0$.
    For each $\alpha \in (0,1)$, let 
    $$
    n(\alpha) := \min\left\{n :  \forall n' \ge n \text{, } \left|\frac{f_0\left(n',\floor{\alpha n'},\floor{c_2{n'-1\choose 2}}\right)}{{n'-1\choose 2}}-C\left(\alpha, \frac{c_2}{2}\right)\right| \le \varepsilon\right\}
    $$
    The value $n(\alpha)$ exists due to the convergence of the function $\frac{f_0\left(n,\floor{\alpha n},\floor{c_2{n-1\choose 2}}\right)}{{n-1\choose 2}}$ proven in Lemma \ref{lem:discrete-continuous-f0-g}. Let $n_0 := \sup_{\alpha\in (0,1)}\left\{n(\alpha)\right\}$. This is $O(\poly(\frac{1}{\varepsilon}))$ due to Lemma~\ref{lem:discrete-continuous-f0-g}.
    
%    \st{Fix a $\varepsilon > 0$.}
    We claim that for all $n\ge n_0$,
    $$
    \left|\frac{f^*(n,\floor{c_2 {n-1\choose 2}})}{{n-1\choose 2}} - c_1^*(c_2) \right| \leq \varepsilon.
    $$
    Indeed, let $\alpha^* := \arg \max_{\alpha \in (0,1)} \{C\left(\alpha,\frac{c_2}{2}\}\right)$ and let $k = \floor{\alpha^* n}$ (for $n \geq n_0$). Then 
    $$
    c_1^*(c_2) - \frac{f_0(n,\floor{\alpha^* n},\floor{c_2 {n-1\choose 2}})}{{n-1\choose 2}} \le \varepsilon
    $$
    since $c_1^*(c_2) = C\left(\alpha^*,\frac{c_2}{2}\right)$ by definition. In particular, this means that $$\frac{\max_{k}\left\{ f_0\left(n,k,\floor{c_2 {n-1\choose 2}}\right)\right\}}{{n-1 \choose 2}} \geq c_1^*(c_2) - \varepsilon.$$
    
    Now, for any $k'$, let $\alpha := \frac{k'}{n}$. Then 
    $$
    \frac{f_0(n,\floor{\alpha n}, \floor{c_2 {n-1\choose 2}})}{{n-1\choose 2}} \le C\left(\alpha, \frac{c_2}{2}\right) + \varepsilon \le C\left(\alpha^*,\frac{c_2}{2}\right) + \varepsilon = c_1^*(c_2) + \varepsilon. 
    $$ 
    Therefore,
    $$
    \frac{f_0(n,\floor{\alpha n}, \floor{c_2 {n-1\choose 2}})}{{n-1\choose 2}} - c_1^*(c_2) \le \varepsilon.  
    $$
    Thus we conclude that $\frac{\max_{k}\left\{ f_0\left(n,k,\floor{c_2 {n-1\choose 2}}\right)\right\}}{{n-1\choose 2}}$ cannot be smaller than $c_1^*(c_2) - \varepsilon$ and cannot be larger than $c_1^*(c_2) + \varepsilon$.
   
   The proof of the limit for $\frac{g^*(n,\floor{c_2 {n-1\choose 2}})}{{n-1\choose 2}}$ is analogous. \end{proof}

\subsection{Proof of Lemma~\ref{lem:o-n-to-1-plus-epsilon}}

\begin{proof}
%    First we prove that for any $\alpha\in (0,1)$, $\lambda > 0$ and $\varepsilon > 0$, there exists an $n_0 =  poly(\frac{1}{\alpha}, \frac{1}{\lambda}, \frac{1}{\varepsilon})$ such that for all $n\ge n_0$,
%    $$
%    C\left(\alpha, \frac{c_2}{2}\right) {n-1\choose 2} - \lambda n^{1+\varepsilon} \le f_0\left(n,\lfloor\alpha n\rfloor,\left\lfloor c_2 {n-1\choose 2}\right\rfloor\right).
%    $$
    Let $\alpha^* := \arg\max_{\alpha \in (0,1)}\left\{C\left(\alpha, \frac{c_2}{2}\right)\right\}$. 
    From the definition of $f^*(n,\left\lfloor c_2{n-1\choose 2}\right\rfloor)$, we know that
    $$
    f_0\left(n,\lfloor \alpha^* n\rfloor, \left\lfloor c_2{n-1\choose 2}\right\rfloor\right) \le f^*\left(n,\left\lfloor c_2{n-1\choose 2}\right\rfloor\right).
    $$
    Therefore, we are going to prove that for any $0<c_2<1$, $\lambda >0$ 
    and  $\varepsilon > 0$, there exists an $n_0$ such that for all $n\ge n_0$,
    $$
     c_1^*(c_2){n-1\choose 2} - \lambda 
     n^{1+\varepsilon} \le  f_0\left(n,\lfloor \alpha^* n\rfloor, \left\lfloor c_2{n-1\choose 2}\right\rfloor\right),
    $$
    which proves the lemma. 
    
    There are $3$ cases. We prove the lemma for all cases.
    
    (\emph{Case 1.}) If $\frac{\alpha^*}{2} < \frac{c_2}{2} <\alpha^*\left(1-\frac{\alpha^*}{2}\right)$ then $c_1^*(c_2) = \frac{\alpha^*}{1-\alpha^*}\left(\frac{c_2-(\alpha^*)^2}{2}\right)$. Also, for sufficiently large $n$, condition $[C1]$ holds with $k = \an$ and $\dmax = \cn$, and thus    
%    ${\lfloor \alpha^* n\rfloor-1 \choose 2} < c_2 {n-1\choose2} \le {\lfloor \alpha^* n\rfloor-1 \choose 2} + (n-\lfloor\alpha^* n\rfloor)(\lfloor\alpha^* n\rfloor -1)$ and thus 
    $$
    f_0\left(n,\an, \cn\right) = 
    $$
    $$
    =\frac{1}{n-\an}\left(\left\lfloor\frac{\an \left(\cn-{\an -1 \choose 2}\right)}{2}\right\rfloor\right. +$$
    $$
   \left. +2\left( \an \left(\cn-{\an -1 \choose 2}\right) (\mod 2)\right)\right).
    $$
    Therefore,
    $$
     c_1^*(c_2){n-1\choose 2} - f_0\left(n,\an,\cn\right) \le
     $$
     $$
     \le \frac{\alpha^*}{1-\alpha^*}\frac{c_2-(\alpha^*)^2}{2} \frac{n^2}{2} -  \frac{1}{n-\alpha^* n} \left(\frac{(\alpha^* n -1)(c_2{n-1\choose 2}-1-\frac{(\alpha^* n)^2}{2})}{2}-1\right) =
    $$
    $$
    = n^2 \left(\frac{\alpha^*}{1-\alpha^*}\frac{c_2-(\alpha^*)^2}{4} - \frac{\alpha^*-\frac{1}{n}}{1-\alpha^*}\left(\frac{c_2(1-\frac{3}{n} + \frac{2}{n^2})-\frac{1}{n^2}-(\alpha^*)^2}{4}-\frac{1}{n^2}\right)\right) =
    $$
    $$
    \frac{1}{4(\alpha^*-1)}\left(\left((\alpha^*)^2-3\alpha^*-c_2\right)n+2\alpha^*c_2-5\alpha^*+3c_2-\frac{2c_2-5}{n}\right).
    $$
    Since this expression is $O(n)$, for any $\lambda >0$ and $\varepsilon>0$ it will be smaller than $\lambda n^{1+\varepsilon}$ if $n$ is sufficiently large.

    (\emph{Case 2.}) If $\frac{c_2}{2} > \alpha^*\left(1-\frac{\alpha^*}{2}\right)$ then $c_1^*(c_2) = \frac{2\alpha^*}{1-\alpha^*}(c_2-(\alpha^*)^2) - 3(\alpha^*)^2$. Also, for sufficiently large $n$, condition $[C2]$ holds with $k= \an$ and $\dmax = \cn$, and thus
    $$
    f_0\left(n,\an\cn\right) = {\an \choose 2} +
    $$
    $$
    +\frac{1}{n-\an}\left(2\an (\cn-{\an -1 \choose 2} - (n-\an)(\an -1))\right).
    $$
    Therefore,
     $$
     c_1^*(c_2){n-1\choose 2} - f_0\left(n,\an,\cn\right) \le
     $$
    $$
    \le \left(\frac{2\alpha^*}{1-\alpha^*}(c_2-(\alpha^*)^2) - 3(\alpha^*)^2\right)\frac{n^2}{2} -\frac{(\alpha^*n-2)^2}{2}- 
    $$
    $$
    -\frac{2\alpha^*-\frac{2}{n}}{1-\alpha^*}\left(c_2\frac{(n-1)(n-2)}{2}-\frac{(\alpha^*)^2n^2}{2}-(n-\alpha^* n)\alpha^*n\right)\le
    $$
    $$
    \le n^2\left(\frac{\alpha^*}{1-\alpha^*}\left((c_2-(\alpha^*)^2)\right) - \frac{3(\alpha^*)^2}{2}-\frac{(\alpha^*)^2-\frac{3\alpha^*}{n}+\frac{2}{n^2}}{2}-\right.
    $$
    $$
   -\left.\frac{2\alpha^*-\frac{2}{n}}{1-\alpha^*}\left(\frac{c_2\left(1-\frac{3}{n}+\frac{2}{n^2}\right)}{2}-\frac{(\alpha^*)^2}{2}-\alpha^*(1-\alpha^*)\right)\right) =
    $$
    $$
    = \frac{1}{2(\alpha^*-1)}\left((\alpha^*)^2n-6\alpha^*c_2n+\alpha^*n-2c_2n+4\alpha^*c_2-2\alpha^*+6c_2+2-\frac{4c_2}{n}\right)
    $$
  Again, this expression is $O(n)$, thus, for any $\lambda > 0$ and $\varepsilon>0$ it will be smaller than $\lambda n^{1+\varepsilon}$ if $n$ is sufficiently large.

(\emph{Case 3.}) If $\frac{c_2}{2} = \alpha^*\left(1-\frac{\alpha^*}{2}\right)$ then $c_1^*(c_2) = \frac{\alpha^*}{1-\alpha^*}\left(\frac{c_2-(\alpha^*)^2}{2}\right) = (\alpha^*)^2$. Observe the following: if $c_2 = \alpha^*(2-\alpha^*)$, then also $\frac{2\alpha^*}{1-\alpha^*}(c_2-(\alpha^*)^2) - 3(\alpha^*)^2 = (\alpha^*)^2$. Therefore, it does not matter if condition [C1] or [C2] holds with $k = \left\lfloor\alpha^*n\right\rfloor$ and $d_{max} = \left\lfloor c_2 {n-1 \choose 2}\right\rfloor$, in both cases $ c_1^*(c_2){n-1\choose 2} - f_0\left(n,\an,\cn\right)$ will be an $O(n)$ function (as shown in Case 1 and 2), and thus, for any $\lambda >0$ ans $\varepsilon > 0$, it will be smaller than $\lambda n^{1+\varepsilon}$ if $n$ is sufficiently large.

\end{proof}

\section{Omitted Proofs from Section \ref{sec:asymp}}\label{sec:appendix-proofs}

\begin{proof}[Proof of Lemma \ref{lem:asymp-t-prob}]
Denote $X_A := |W \cap A|$ and $X_B := |W \cap B|$. Because $W = (W \cap A) \sqcup (W \cap B)$ and $|W| = t$, we have the following equality between events:  

$$ \left\{\left||W \cap A|-|W \cap B|\right|> \delta t \right\} = \ 
\left\{\frac{\left|X_A-(t - X_A)\right|}{t}> \delta \right\} = \ 
\left\{\left|\frac{X_A}{t} - \frac{1}{2}\right|> \frac{\delta}{2} \right\}.
$$

Thus we want to show that $\pr{\left|\frac{X_A}{t} - \frac{1}{2}\right|> \frac{\delta}{2}} < \epsilon$. Let $\hyp{n}{s}{k}$ denote the hypergeometric distribution with $n$ total items, $s$ success items, and $k$ draws ($k \leq n$). Let $m := \frac{n}{2}$ for convenience. We can view $W$ as being obtained through $t$ uniformly random draws without replacement from a set of size $|V| = 2m$. By viewing $A$ as the set of ``success'' cases of size $m$, it follows that the random variable $X_A$ has a hypergeometric distribution; in particular, $X_A \sim \hyp{2m}{m}{t}$. 

\vspace{0.5em}
We can now apply Chebyshev's inequality with the hypergeometric distribution to conclude our bound, recalling that for $X \sim \hyp{n}{s}{k}$, $\Ex[X] = k \cdot \frac{s}{n}$ and $\Var[X] = k \cdot \frac{s}{n} \cdot \frac{n-s}{n} \cdot \frac{n-k}{n-1}$: 

\begin{align*}
\pr{\left|\frac{X_A}{t} - \frac{1}{2}\right|> \frac{\delta}{2}} 
&= \pr{\left|\frac{X_A}{t} - \Ex\left[\frac{X_a}{t}\right]\right|> \frac{\delta}{2}} & \text{ since $X_A \sim \hyp{2m}{m}{t}$} \\ 
&\leq \Var\left[\frac{X_A}{t}\right] / \left(\frac{\delta^2}{4}\right) & \text { by Chebyshev's inequality} \\ 
&= \frac{4 \Var[X_A]}{\delta^2 t^2} \\ 
&< \frac{1}{\delta^2 t} & \text{$\Var[X_A] = \frac{t}{4} \cdot \frac{2m-t}{2m-1} < \frac{t}{4}$} \\
&\leq \epsilon & \text{$\delta \geq (\epsilon \cdot t)^{-\frac{1}{2}}$}
\end{align*}\end{proof}

% {\color{black!50}
% \section{[Remove Later] Flowchart: Results in Section \ref{sec:always=graphic}}
% \begin{itemize}
%     \item Lemma \ref{lem:min-f}: uses Definition \ref{def:f}
%     \item Lemma \ref{lem:g-int-graphic} uses Definition \ref{def:f}, Lemma \ref{lem:min-f} (Part 2), Definition \ref{def:g}
%     \item Lemma \ref{lem:n-k=o(n)}: uses Definition \ref{def:f} (and the ``add-and-HF sub-routine''
%     \item Lemma \ref{lem:linear-bounds}: Lemma \ref{lem:g-int-graphic}, Lemma \ref{lem:n-k=o(n)} (\textit{Check: don't think we directly use Lemma \ref{lem:min-f} (Part 2) here?})
%     \item Theorem \ref{thm:symmetric-bounds}: uses Definition \ref{def:f}, Lemma \ref{lem:linear-bounds}
% \end{itemize}}

\end{document}

%% file: commands.tex
\newcommand{\grparam}{\textsc{3uni-HDS}_{c_1,c_2}}
\newcommand{\grprob}{\textsc{3uni-HDS}}
\newcommand{\gr}{ k\textsc{uni-HDS} }

\newcommand{\dcc}{\mathcal{D}_{c_1,c_2}}
\newcommand{\generaldcce}{\mathcal{D}_{c_1,c_2}(\varepsilon)}
\newcommand{\dcce}{\mathcal{D}_{c_1^*,c_2}(\varepsilon)}

\newcommand{\dmin}{d_{min}}
\newcommand{\dmax}{d_{max}}
\newcommand{\dint}{d_{int}}
\newcommand{\poly}{\text{poly}}

\newcommand{\an}{\lfloor \alpha^* n \rfloor}
\newcommand{\cn}{\left\lfloor c_2 {n-1\choose 2} \right\rfloor}

\newcommand{\EH}{E_{high}}
\newcommand{\EL}{E_{low}}

\newcommand{\hyp}[3]{\textsc{HGeom}({#1}, {#2}, {#3})}

\newcommand{\pmodtight}[2]{{#1} \!\!\!\! \pmod{#2}}